\documentclass[11pt,reqno,fleqn]{amsart}

\spaceskip=1.1\fontdimen2\font plus 1.0\fontdimen3\font minus 0.5\fontdimen4{amsmath}
\usepackage{amsfonts}
\usepackage{amssymb}
\usepackage{amsgen}
\usepackage{enumitem}
\usepackage{graphicx}
\usepackage{hyperref}
\usepackage{mathrsfs}
\usepackage{cite}

\setlist[enumerate]{label=\arabic*.,left=1em}

\newdimen\smmathindent
\usepackage{latexsym}

\usepackage{textcomp}
\usepackage[sb]{libertine}
\usepackage[T1]{fontenc}
\usepackage[varqu,varl]{zi4}
\usepackage[libertine,bigdelims,vvarbb]{newtxmath} 
\usepackage[cal=boondoxo,frak=euler]{mathalpha} 
\usepackage[supsfam=LibertinusSerif-Sup,%
	supscaled=1.2,%
	raised=-.13em]{superiors}

\newcommand{\term}[1]{\emph{#1}}

\def\nnsection#1\par{\medbreak%
	\penalty-5000%
	\\{\large\bf #1}\medbreak%
	}

\newcommand{\customhead}[1]{\medbreak\\ \textsc{\underline{\smash{#1\!}\phantom{$,$}}}\medbreak}

\newcommand{\q}{\quad}
\newcommand{\qq}{\q\q}

\renewcommand{\\}{\noindent\ignorespaces}

\newcommand{\Z}{\mathbb{Z}}
\newcommand{\Q}{\mathbb{Q}}
\newcommand{\R}{\mathbb{R}}
\newcommand{\C}{\mathbb{C}}

\newcommand{\ndn}{\noindent\hskip\the\mathindent}
\renewcommand{\[}{\medbreak\begingroup\ndn$\displaystyle\begingroup}
\renewcommand{\]}{\endgroup$\endgroup\medbreak}

\def\ul#1\\{\ndn#1\par}
\def\uul#1\luu{\ndn#1\par}

\def\nl#1 #2\\{\\\hbox to \the\mathindent{\mathsurround=0pt#1\hfil}#2\par}

\def\cl#1 #2\\{\medbreak\\\hbox to \the\mathindent{\mathsurround=0pt#1\hfil}#2\medbreak}
\def\ccl#1 #2\lcc{\medbreak\\\hbox to \the\mathindent{\mathsurround=0pt#1\hfil}#2\medbreak}

\newcommand{\Aut}{\mathop{\mathsf{Aut}}}

\newcommand{\idrep}{\mathbf{I}}
\newcommand{\subrep}{\mathbf{S}}

\newcommand{\B}{\mathcal{B}}
\renewcommand{\C}{\mathcal{C}}
\newcommand{\F}{\mathcal{F}}
\newcommand{\G}{\mathcal{G}}

\renewcommand{\a}{\alpha}
\renewcommand{\b}{\beta}

\renewcommand{\l}{\lambda}
\newcommand{\s}{\sigma}
\newcommand{\Om}{\Omega}

\newcommand{\fK}{\mathfrak{K}}
\newcommand{\fL}{\mathfrak{L}}

\newcommand{\bk}{\mathbb{k}}
\newcommand{\Bk}{\B_{\bk}}
\newcommand{\BQ}{\B_\Q}

\renewcommand{\phi}{\varphi}

\newcommand{\setin}{\subseteq}
\newcommand{\contains}{\supseteq}

\newcommand{\Lax}{\mathsf{Lax}}

\newcommand{\ER}{\textbf{ER}}

\newcommand{\ssum}{\vee}
\newcommand{\sprod}{\wedge}

\newcommand{\cO}{\mathcal{O}}
\newcommand{\cS}{\mathcal{S}}
\newcommand{\cT}{\mathcal{T}}
\newcommand{\cs}{\mathcal{s}}

\newcommand{\rank}{\mathop{\rm rank}}
\newcommand{\im}{\mathop{\rm im}}
\renewcommand{\ker}{\mathop{\rm ker}}
\newcommand{\id}{\mathord{\rm id}}

\newcommand{\Sym}{\mathop{\rm Sym}}

\newcommand{\rac}{\ast}
\newcommand{\lac}{\ast}

\newcommand{\scr}[1]{\mathbin{\mathscr{#1}}}

\renewcommand{\R}{\scr{R}}
\renewcommand{\L}{\scr{L}}
\newcommand{\J}{\scr{J}}
\renewcommand{\H}{\scr{H}}

\newcommand{\ger}{\ge_{\scr{R}}}
\newcommand{\gel}{\ge_{\scr{L}}}
\newcommand{\gej}{\ge_{\scr{J}}}
\newcommand{\geh}{\ge_{\scr{H}}}

\newcommand{\emax}{\equiv_{\mathord{\max}}}

\renewcommand{\B}{\mathcal{B}}

\newcommand{\Thmname}{Theorem}
\newcommand{\Propname}{Proposition}
\newcommand{\Lemmaname}{Lemma}
\newcommand{\Definitionname}{Definition}
\newtheorem{Thm}{\Thmname}[section]
\newtheorem{Prop}[Thm]{\Propname}
\newtheorem{Lemma}[Thm]{\Lemmaname}
{\theoremstyle{definition}
}
{\theoremstyle{remark}
}
\newtheorem{Cor}[Thm]{Corollary}
{\theoremstyle{remark}
}

\newtheorem*{Lemma*}{Lemma}

\numberwithin{equation}{section}

\title{A Burnside ring for monoids}

\author{Jeremy Weissmann}
\email{jweissmann@gradcenter.cuny.edu}

\begin{document}
\begin{abstract}
In this note, we define the Burnside ring of a monoid, generalizing the construction for groups.  After giving foundational definitions, we describe the monoid-theoretic correlates of orbits and their automorphisms, then prove a structure theorem for a broad class of monoids that allows us to write the Burnside ring as a direct product of Burnside rings of groups. Finally, we define a monoid-theoretic correlate of the table of marks, and show that the Burnside algebra over $\Q$ is semisimple.
\end{abstract}

\maketitle
\tableofcontents

\renewcommand{\baselinestretch}{1.1}
\setlength{\parskip}{4pt}

\setcounter{section}{-1}
\section{Introduction}
	
	The Burnside ring of a finite group \cite{burnside,solomon} is an algebraic representation of all the ways that group can act on sets, and as such, it is a fundamental tool in representation theory.  In this paper, we define the Burnside ring of a finite monoid, generalizing the group-theoretic construction.  We also prove several monoid-theoretic correlates of standard results about the Burnside ring.

	As a preview, let us outline the construction of the Burnside ring for groups:  First, one defines what it means for a group to \emph{act} on a set; a set acted on by a group~$G$ is called a \emph{$G$-set}. One defines a notion of isomorphism between $G$-sets, and the elements of the Burnside ring are isomorphism classes of $G$-sets. Addition is given by disjoint union; multiplication by Cartesian product.  It is easily verified that these operations yield a well-defined semiring. Finally, one uses a Grothendieck-style construction to create additive inverses, and embed the semiring in a ring.
	
	All of this can be translated into monoid-theoretic terms:  A set acted on by a monoid~$M$ is called an \emph{$M$-set}, and so forth.  This would lead us to a construction similar to what we will call the \emph{weak Burnside ring} of a monoid in Section~2.  As this name suggests, however, the weak Burnside ring does not quite capture the structure of a monoid in the way one would hope.

	Let us briefly explain the difficulty:  If a group element $g$ takes $x$ to $y$ via the group action, then $g^{-1}$ takes $y$ to $x$. Consequently, a $G$-set can be decomposed into \emph{orbits}, subsets which are \emph{invariant} (closed under the group action) and \emph{transitive} (connected via the group action). By contrast, since monoid elements are not necessarily invertible, it may be possible to go from~$x$ to~$y$ via the action but not back, so that the ``orbit'' containing~$x$ must contain $y$ by invariance, but then transitivity fails.
	
	As we will see, maximally connected sets, which we will call \emph{strong orbits}, are the appropriate correlate of orbits in the group-theoretic setting.  However, strong orbits are not necessarily invariant, so we cannot decompose $M$-sets into strong orbits.  We are thus led to consider a more general notion of actions \emph{by partial functions}. Once we have generalized in this way, we will be able to impose certain relations on the weak Burnside ring so that such decompositions are possible.  Thus the Burnside ring for monoids arises as a quotient of the weak Burnside ring.

\section{Foundations}

	In this section, we present the concepts, definitions, and notations which are foundational to the construction of the Burnside ring of a monoid.  At the heart of the construction is the notion of a monoid acting on a set by partial functions.

\subsection{Partial functions and definedness}

	We begin by introducing some terminology which is hopefully familiar to the reader.

	A \term{partial function} $f : X \to Y$ is a relation on sets $X$ and $Y$ which is a function when restricted to its domain.  If $x \in X$ is not in the domain of $f$ we say that $f(x)$ is \term{undefined}.
		
	We stipulate that \emph{an expression is undefined if it depends functionally on an undefined expression}. Conversely, \emph{every subexpression of a defined expression is defined}. This is not a philosophical position: it will have real significance when we define subquotients below in Section~1.4.

	We have a few more conventions that should be understood when dealing with expressions that may or may not be defined: When we write, for instance, that $f(x) \in Y$, we mean that this set membership holds where $f(x)$ is defined.  More subtly, if we write $f(x) = f(x')$, this is intended to mean that both expressions are defined and equal, or both are undefined.  And so forth.

	Partial functions are only used in the definition of an action, given in the next section.  All other functions, maps, operations, etc, are assumed total, with the exception of the discussion in Section~1.7, which is self-contained.

\subsection{Monoid actions and $M$-sets}

	A \term{semigroup} is a set closed under an associative binary operation; a \term{monoid} is a semigroup with an \term{identity} element; and a \term{group} is a monoid which contains an \term{inverse} for each element. We use multiplicative notation for these structures throughout the paper.
		
	We focus mainly on monoids, and consider groups as a special case.  We denote an arbitrary monoid by~$M$, and in the case that $M$~is a group, we let $M = G$ for clarity.

	An \term{action by partial functions} of $M$ on a set $X$ associates each element of~$M$ with a partial function from $X$ to itself, in a way that respects the structure of~$M$. More precisely, $\rac$ denotes a \term{right action by partial functions} if for all $x \in X$ and $m,n \in M$, the following axioms hold:
	
	\medbreak
	\nl $\q\bullet$ \hbox to 13.5em{$x \rac m \in X$ \ where defined\hfil} {\bf closure}\\
	\nl $\q\bullet$ \hbox to 13.5em{$(x \rac m) \rac n \,=\,
		x \rac mn$\hfil} {\bf associativity}\\
	\nl $\q\bullet$ \hbox to 13.5em{$x \rac 1 = x$\hfil} {\bf action of the identity.}\\
	\medbreak

\\ As discussed above, associativity should be understood as saying that both sides of the equation are defined and equal, or both are undefined.  Note that associativity implies that if $x \rac mn$ is defined, then $x \rac m$ is defined as well.  The action of the identity means that $x \rac 1$ is always defined and equal to~$x$.
	
	Left actions are defined dually.  Most of the actions under discussion are right actions, but an important left action will be introduced in Section~3.4.

	If an action as defined above exists, we say that \term{$M$ acts on $X$ by partial functions}, and that $X$ is a \term{partial $M$-set}. However, because we have to use these terms so frequently, we drop the word `partial' and simply say that $M$~acts on $X$ and that $X$ is an $M$-set. If $x \rac m$ is defined for all $x \in X$ and $m \in M$, we say that the action is \term{total} and that $X$ is a \term{total $M$-set}.
	
	Group actions are total.  Indeed, if $M = G$ is a group, then:
	
	\[
	x \,=\, x \rac 1 \,=\, x \rac gg^{-1}
	\]
	
\\ is defined, and hence $x \rac g$ is defined for all $x \in X$ and $g \in G$. We will carry out many calculations like this in the sequel, and it is worth drawing out the principle: {\em If we begin a calculation with a defined expression, then all subsequent expressions are defined as well.}

	Finally, we often write $xm$ instead of $x \rac m$.

\subsection{Invariant subsets}

	Let $X$ be an $M$-set, and let $A$ be a subset of $X$. We say that $A$ is \term{invariant} if it is closed under the action of~$M$. In symbols, this means that for all $a \in A$ and $m \in M$ we have $am \in A$ where defined.  If $A$ is invariant, then $A$ is an $M$-set under the action of $M$ restricted to~$A$, because the other axioms are satisfied automatically.
	
	Letting $AM = \{ am : a \in A,\, m \in M,\,am\ {\rm defined} \}$, the invariance of~$A$ may be expressed by~$AM \setin A$. Note that for any $A \setin X$ we have:
	
	\[
	(AM)M = A(MM) \setin AM
	\]
	
\\ so all subsets of the form $AM$ are invariant.

	By convention, we write $aM$ instead of $\{a\}M$.

\subsection{Restrictions and subquotients}

	Let $M=G$ be a group, and let $A$ be an invariant subset of a $G$-set $X$. If $xg \in A$, then $x = x1 = xgg^{-1} \in A$ by invariance. By the contrapositive, $x \in X \setminus A$ implies $xg \in X \setminus A$; in other words, the complement $X \setminus A$ is invariant as well.  This property, alluded to in the introduction, allows one to decompose $G$-sets into orbits.
	
	This property does not hold for monoids in general.  For example, let $M = \{0,\pm1\}$ under multiplication, and let $M$ act on itself by right multiplication.  Then $\{0\}$ is an invariant subset of~$M$, but its complement $\{\pm1\}$ is not, since for example $1 \rac 0 = 0 \notin \{\pm1\}$. However, if we allow $1 \rac 0$ and $-1 \rac 0$ to be undefined (since the value these expressions should have is not in $\{\pm1\}$), it is easy to verify that $\{\pm 1\}$ is a partial $M$-set.  As we will see, this phenomenon holds generally: the complement of an invariant subset is a partial $M$-set.
	
	To be precise:  Suppose $S$ is any subset of an $M$-set $X$. Let $s \cdot m = sm$ if $sm$ is defined and $sm \in S$, and let $s \cdot m$ be undefined otherwise.  We call $\cdot$ the \term{restriction to $S$} of the action on~$X$. It is straightforward to check that $\cdot$ defines an action if and only if:
	
	\[
	smn \in S \q{\rm implies}\q sm \in S
	\]
	
\\ or in terms of sets, that $SM \setminus S$ is invariant. In this case, we say that $S$ is a \term{subquotient} of~$X$.  (The reason for the term \term{subquotient} will be explained in Section~1.7.)

	Intuitively, $S$ is a subquotient if the action on~$X$ cannot take elements out of~$S$ and back in again. If $S$ is invariant, then $S$ is a subquotient, since the action never takes elements out of an invariant set.  But this is equivalent to saying that the action never takes elements into the complement of an invariant set, so $X \setminus S$ is a subquotient as well, which we denote by~$X/S$.
	
	The fact that both invariant sets and their complements are subquotients is what will allow us in Section~2 to decompose $M$-sets analogously to the case of groups.

\subsection{Morphisms}

	A \term{morphism of $M$-sets} is a map which preserves the action of $M$. In symbols, we say that a map $f : X \to Y$ of $M$-sets is a \term{morphism} if $f(xm) = f(x)m$ for all $x \in X$ and $m \in M$. As usual, this means that both expressions are defined and equal, or both are undefined.  The definedness condition is equivalent to showing that $xm$ is defined if and only if $f(x)m$ is.
	
	To illustrate this definition with an important example, let $S$ be a subset of $X$, and define the inclusion map $\iota : S \to X$ by $\iota(s) = s$. If $S$ is invariant then $\iota(sm) = sm = \iota(s)m$ for all $m \in M$, so $\iota$ is a morphism. However, if $S$ is a subquotient then $\iota$ is a morphism if and only if $\iota(s \cdot m) = \iota(s)m = sm$, if and only if $s \cdot m = sm$ for all $m \in M$. In other words, inclusion of a subquotient is a morphism if and only if the subquotient is invariant. Inclusion of a subquotient is an example of what we will call a \term{lax morphism} in Section~5.

	We define the \term{image} of $f$ by:
	
	\[
	\im(f) \,=\, f(X) \,=\, \{ f(x) : x \in X \}.
	\]

\\ The morphism property tells us that $f(X)M = f(XM) \setin f(X)$, so $\im(f)$ is invariant.

	If $f$ is a bijection, we say that $f$ is an \term{isomorphism} and that $X$ and $Y$ are \term{isomorphic}, and write $X \cong Y$. It is clear that $\cong$ is an equivalence relation.  If $f$ is an isomorphism, it is straightforward to check that its inverse $f^{-1}$ is a morphism.

\subsection{Congruence and quotients}
	
	Let $X$ be an $M$-set. We call an equivalence relation $\equiv$ on $X$ a \term{congruence} if it is preserved by the action of $M$, meaning that $x \equiv x'$ implies $xm \equiv x'm$ for all $x,x' \in X$ and $m \in M$. As usual, this means that both expressions are defined and equivalent, or both are undefined. If $\equiv$~is a congruence, we speak of \term{congruence classes} instead of equivalence classes.
			
	Let $\equiv$ be a congruence on $X$, let $X/{\equiv}$ be the set of congruence classes, and let $\overline x$ denote the congruence class of~$x$. Then $M$ acts on $X/{\equiv}$ by letting $\overline{x}m = \overline{xm}$ if $xm$ is defined, and letting $\overline{x}m$ be undefined otherwise. We call $X/{\equiv}$ a \term{quotient of $X$}. The \term{quotient map} from $X$ to $X/{\equiv}$ given by $x \mapsto \overline{x}$ is a surjective morphism.
	
	Let $f : X \to Y$ be a morphism, and define $\equiv_f$ on $X$ by $x \equiv_f x'$ if and only if $f(x) = f(x')$. It follows that $\equiv_f$ is a congruence, and we have the usual isomorphism $X/{\equiv_f} \cong \im(f)$ via the map $\overline{x} \mapsto f(x)$.

\subsection{Pointed $M$-sets}

	In this section, we discuss an equivalent approach to the foundations, in order to connect the dots between our work and Haesemeyer and Weibel's $K$-theoretic approach in~\cite{weibel}.  The concepts and definitions introduced in this section will not be used elsewhere in this paper.
	
	A \term{pointed set} is a pair $(X,0_X)$ where $X$ is a set and $0_X$ is an element of~$X$ called the \term{distinguished point}. We say that a monoid $M$ \term{acts} on a pointed set $(X,0_X)$ via $\rac$ if for all $x \in X$ and $m,n \in M$, the following axioms hold:
	
	\medbreak
	\nl $\q\bullet$ \hbox to 12em{$x \rac m \in X$\hfil} {\bf closure}\\
	\nl $\q\bullet$ \hbox to 12em{$(x \rac m) \rac n \,=\,
		x \rac mn$\hfil} {\bf associativity}\\
	\nl $\q\bullet$ \hbox to 12em{$x \rac 1 = x$\hfil} {\bf action of the identity}\\
	\nl $\q\bullet$ \hbox to 12em{$0_X \rac m \,=\, 0_X$\hfil} {\bf action on the distinguished point.}\\
	\medbreak

\\ In this case, we say that $\rac$ is a \term{pointed action} and that $(X,0_X)$~is a \term{pointed $M$-set}. Note that we require $x \rac m$ to be defined for all $x \in X$ and $m \in M$. For simplicity we write $xm$ instead of $x \rac m$.

	If we take $xm = 0_X$ as saying that $xm$ is undefined, we see that pointed and partial $M$-sets describe the same concept with different language.  Note that closure no longer holds ``where defined'', but simply for all $x \in X$ and $m \in M$; and similarly, we do not need to stipulate that associativity means that both sides are defined and equal or both are undefined. Finally, observe that the action on the distinguished point corresponds to our stipulation that expressions which depend on undefined expressions are undefined.
	
	Let $\equiv$ be an equivalence relation defined on a pointed $M$-set $X$. If $x \equiv x'$ implies $xm \equiv x'm$ for all $x,x' \in X$ and $m \in M$, we say that $\equiv$ is a \term{right congruence}. The \term{quotient} $X/{\equiv}$ is a pointed $M$-set where the class of $0_X$ is the distinguished point.
	
	We say that $A$ is an \term{invariant} subset of a pointed $M$-set $X$ if $AM \setin A$ and $0_X \in A$. If $A$~is invariant, we may define a congruence $\equiv_A$ by $x \equiv_A y$ if and only if $x = y$ or $x,y \in A$, so that $X/{\equiv}_A$ is the quotient obtained by crushing~$A$ to the distinguished point.  We write this quotient as $X/A$. Recall that in the partial $M$-set approach, what we wrote as $X/A$ was actually just the complement $X \setminus A$ under the restriction of the action.  By contrast, in the pointed approach we really have a quotient.  This explains why we chose the term \term{subquotient}: If $SM \setminus S$ is invariant, then $S$ is a quotient of the sub-pointed-$M$-set~$SM$, namely $SM/(SM \setminus S)$.

	A \term{morphism} of pointed $M$-sets is a map $f : (X,0_X) \to (Y,0_Y)$ such that $f(0_X) = 0_Y$ and $f(xm) = f(x)m$ for all $x \in X$ and $m \in M$. Define the \term{image} by:
	
	\[
	\im(f) \,=\, f(X) \,=\, \{ f(x) : x \in X \}
	\]
	
\\ and the \term{kernel} by:
	
	\[
	\ker(f) \,=\, f^{-1}(0_Y) \,=\, \{ x \in X : f(x) = 0_Y \}.
	\]
	
\\ The image and kernel are invariant subsets of $Y$ and $X$, respectively.

	In the language of partial $M$-sets, $f(x) = 0_Y$ means that $f(x)$~is undefined.  This is something new, because above we only considered morphisms of partial $M$-sets which were total functions.  Define $x \equiv_f x'$ if and only if $f(x) = f(x')$. Then $X/{\equiv}_f$ is a pointed $M$-set with distinguished point $\ker(f)$, and $X/{\equiv}_f$ is isomorphic to $\im(f)$ as pointed $M$-sets.
	
	Recall that in the partial $M$-set approach, we also had $X/{\equiv_f} \cong \im(f)$. However, this only held for total morphisms.  If $f$ can be partial, then $X/{\equiv_f}$ contains the congruence class $K = \{ x \in X : f(x)\ {\rm undefined}\}$ which has no correlate in~$\im(f)$.  One might try to address this by restricting $\equiv_f$ to the domain of~$f$; however, the domain of~$f$ is not necessarily invariant.  Instead we must use the fact that $\{K\}$ is an invariant subset of $X/{\equiv_f}$ so that $(X/{\equiv_f})/\{K\}$ is a subquotient isomorphic to~$\im(f)$. Not particularly nice.
	
	So we see that, while the partial and pointed approaches are equivalent, the pointed approach has certain advantages: it is a little more homogeneous and allows us to more easily discuss what we would call partially-defined morphisms.  But like any approach, it has disadvantages as well:  In particular, the definitions of sums and products of $M$-sets, to be introduced shortly, become more cumbersome, which in turn makes the whole construction of the Burnside ring less pleasant.  Since this is after all a paper on the Burnside ring, and since we do not need the generality of partially-defined morphisms, we will follow the partial $M$-set approach exclusively in what follows.

\section{The Burnside rings of a monoid}

	In this section, we define operations on $M$-sets which will become addition and multiplication in the Burnside ring.  We define \term{weak} and \term{strong orbits}, and give some of their properties with respect to these operations.
	
	The elements of the Burnside ring are equivalence classes of $M$-sets. As discussed in the Introduction, we define two Burnside rings: the \term{weak Burnside ring} $\B_w$, which is freely generated as an additive group by the weak orbit classes, and the \term{strong Burnside ring} $\B$, which is freely generated as an additive group by the strong orbit classes.
	
	We explain why these rings are the same in the case that $M$ is a group, and also suggest why $\B$ is likely to be the more useful construction.

\subsection{Operations on $M$-sets}

\subsubsection{Sums}

	The \term{sum} of two $M$-sets is their disjoint union.  Indeed, if $X$ and~$Y$ are $M$-sets, then $M$~acts on the disjoint union of $X$ and~$Y$, using the appropriate action depending on which set the element comes from.
	
	We denote the sum of $X$ and $Y$ by $X \ssum Y$. (Rather than giving any specific set-theoretic apparatus to deal with disjoint union, we simply assume that $X$ and~$Y$ are disjoint whenever we write~$X \ssum Y$.)
	
	What we have defined is an \term{external sum}. We can look at things the other way as well: Suppose $X$ is an $M$-set, and that a subset $S$ and its complement~$X \setminus S$ are invariant.  Then $X$ is the \term{internal sum} of $S$ and $X \setminus S$, and in this case we also write $X = S \ssum (X \setminus S)$, since it is clear that the external and internal constructions are isomorphic.
	
	If $f : X \to X'$ and $g : Y \to Y'$ are morphisms, then the \term{sum of morphisms} $f \ssum g$ given by:
	
	\[
	(f \ssum g)(s) = \
		\begin{cases}
			\,f(s) & \hbox{if $s \in X$}\\
			\,g(s) & \hbox{if $s \in Y$}
		\end{cases}
	\]

\\ is a morphism from $X \ssum Y$ to $X' \ssum Y'$, and the inclusion maps:

	\medbreak
	\ul \hbox to 8.0em{$\iota_X : X \to X \ssum Y$,\hfil} $\iota_X(x) = x$\\
	\ul \hbox to 8.0em{$\iota_Y : Y \to X \ssum Y$,\hfil} $\iota_Y(y) = y$\\
	\medbreak

\\ are injective morphisms.  (Recall from Section~1 that inclusion of an invariant subset is always a morphism.)

\subsubsection{Products}

	The \term{product} of two $M$-sets is their Cartesian product.  Indeed, if $X$ and $Y$ are $M$-sets, then $M$ acts on $X \times Y$ by letting $(x,y)m = (xm,ym)$ if $xm$ and~$ym$ are both defined, and letting $(x,y)m$ be undefined otherwise.
	
	We denote the product of $X$ and~$Y$ by $X \sprod Y$.

	If $f : X \to X'$ and $g : Y \to Y'$ are morphisms, then the \term{product of morphisms} $f \sprod g$ given by:
	
	\[
	(f \sprod g)((x,y)) \,=\, (f(x),g(y))
	\]

\\ is a morphism from $X \sprod Y$ to $X' \sprod Y'$. However, the projection maps:
	
	\medbreak
	\ul \hbox to 8.0em{$\pi_X : X \sprod Y \to X$,\hfil} $\pi_X((x,y)) \,=\, x$\\
	\ul \hbox to 8.0em{$\pi_X : X \sprod Y \to Y$,\hfil} $\pi_X((x,y)) \,=\, y$\\
	\medbreak

\\ are not morphisms: For example, if $xm$ is defined but $ym$ is undefined, then $\pi_X((x,y)m)$ is not defined, but $\pi_X((x,y))m = xm$ is.  The projection maps are again examples of what we will call \term{lax morphisms} in Section~5.
	
\subsection{Equivalence classes of $M$-sets}

	Interesting properties of $\ssum$ and $\sprod$ only hold up to isomorphism.  For example, $X \sprod Y$ does not equal $Y \sprod X$ in general, but they are isomorphic via the map $(x,y) \mapsto (y,x)$.
	
	For this reason, we shift our consideration to \term{isomorphism classes} of $M$-sets, or \term{classes} for short.  We write $[X]$ for the class of $X$, and define:
	
	\[
	[X]+[Y] \,=\, [X \ssum Y] \qq{\rm and}\qq [X]\cdot[Y] \,=\, [X \sprod Y]. 
	\]
	
\\ These definitions are independent of representatives, for if $f : X \to X'$ and $g : Y \to Y'$ are isomorphisms, then so are:

	\[
	f \ssum g \,:\, X \ssum Y \,\to\, X' \ssum Y' \qq{\rm and}\qq
	f \sprod g \,:\, X \sprod Y \,\to\, X' \sprod Y' .
	\]

\\ It is similarly straightforward to show that $+$ and $\cdot$ are associative and commutative, and that $\cdot$~distributes over~$+$.

	The empty set $\emptyset$ is an $M$-set under the empty action, and it is clear that ${\bf 0} = [\emptyset]$ is the additive identity.
	
	There is a unique total action on a singleton set $\{\s\}$, namely $\s m = \s$ for all $m \in M$. It is easy to check that $\{\s\} \sprod X \cong X$ via the isomorphism $(\s,x) \mapsto x$, hence ${\bf 1} = [\{\s\}]$ is the multiplicative identity.
	
	In the next section, we will show how to decompose an $M$-set into a sum of $M$-sets called weak orbits.  As an infinite $M$-set may necessitate an ``infinite sum'' of weak orbits, we restrict our attention in the sequel to finite $M$-sets.
	
	Therefore, letting $\C$~denote the collection of all equivalence classes of finite $M$-sets, we have shown that $(\C,+)$ is a commutative monoid and that $(\C,+,\cdot)$ is a commutative semiring.

\subsection{Orbits}

	In this section, we discuss $M$-sets called \term{weak} and \term{strong orbits}, terminology introduced by Steinberg in~\cite[Section 2.1]{transformationmonoids}. They are of importance because the classes of these $M$-sets freely generate the additive groups of the weak and strong Burnside rings, respectively.
	
\subsubsection{Weak orbits}

	We say that a nonempty $M$-set is \term{decomposable} if it is isomorphic to the sum of two nonempty~$M$-sets. Since $X$ is assumed finite, we can repeatedly decompose it into smaller $M$-sets until we have written it as a sum of indecomposable $M$-sets. These indecomposable $M$-sets are what we intend to call \term{weak orbits}.

	To define this concept in terms of the action, consider an $M$-set~$X$ and define a relation $\ge$ on the elements of~$X$ by $x \ge y$ if and only if $xm = y$ for some $m \in M$. Write $x \simeq y$ if there exists a sequence $w_0,w_1,\dots,w_n$ of elements of~$X$ such that $w_0 = x$, $w_n = y$, and for all $i$, $0 \le i < n$, either $w_i \ge w_{i+1}$ or $w_{i+1} \ge w_i$. It is clear that $\ge$ is reflexive and transitive, i.e.~a preorder, and that $\simeq$ is the symmetric-transitive closure of $\ge$, and hence an equivalence relation.  We call its equivalence classes the \term{weak orbits} of~$X$. If $X$ comprises a single equivalence class then we say that $X$ is a \term{weak orbit}.
	
	By way of illustration, consider an undirected graph whose vertices are the elements of $X$, and whose edges show which elements of $X$ are related by the action: in other words, $xy$ is an edge if and only if $xm = y$ or $ym = x$ for some $m \in M$. Then $x \simeq y$ means that there exists an undirected walk from~$x$ to~$y$, and the weak orbits of~$X$ are the connected components of this graph.

\subsubsection{Properties of weak orbits}

	Here we collect some simple but important properties of weak orbits.  The auxiliary graph introduced above may be a helpful image to keep in mind when considering the following arguments.

\begin{Prop}\label{woprops}
The following hold:
\begin{enumerate}
\item Weak orbits are invariant.

\item An $M$-set is the (internal) sum of its weak orbits.

\item Weak orbits of isomorphic $M$-sets can be put in a bijective, isomorphic correspondence.

\item The weak orbits of the sum $S \ssum T$ are precisely those of $S$ and $T$, counting multiplicity.

\item Weak orbits are precisely the indecomposable $M$-sets.

\end{enumerate} 
\end{Prop}

\begin{proof}
Let $W$ be a weak orbit.  Let $w \in W$ and $m \in M$ and suppose that $wm$~is defined. Then $w \ge wm$, so $w \simeq wm$, and hence $wm \in W$ as desired.  This proves the first property.
	
	It follows immediately that a finite $M$-set is the internal sum of its weak orbits, since the weak orbits are a finite collection of invariant subsets which partition the $M$-set. This proves the second property.

	If $f$ is an isomorphism, then $x \ge y$ if and only if $f(x) \ge f(y)$. Since $\simeq$ is defined in terms of $\ge$, it follows that the weak orbit of $x$ corresponds isomorphically via~$f$ to the weak orbit of $f(x)$. This proves the third property.

	Consider the weak orbits of the sum $S \ssum T$. Since $S$ and $T$ are invariant, a weak orbit of $S \ssum T$ is contained entirely in~$S$ or in~$T$. It is also clear that $x \simeq y$ holds in $S$ or in $T$ if and only if $x \simeq y$ holds in $S \ssum T$. This shows that the weak orbits of $S \ssum T$ are precisely the weak orbits of $S$ and $T$, proving the fourth property.
	
	Finally, by Property~2, an $M$-set is the sum of its weak orbits, so if a set has more than one weak orbit, it is decomposable. Conversely, suppose that $W$ is a weak orbit and that $W \cong S \ssum T$. By Properties 3 and~4, we must have $S \cong W$ and $T \cong \emptyset$ or vice versa, so $W$ is indecomposable.  This settles the fifth property and concludes the proof.  
\end{proof}

	These properties allow us to prove an important structural theorem.  Recall that $\C$ is the collection of equivalence classes of finite $M$-sets.
	
\begin{Thm}
The weak orbit classes freely generate $(\C,+)$.
\end{Thm}
\begin{proof}	

	Let $\F$ be a free commutative monoid on the weak orbit classes, and let $\phi : \C \to \F$ be the map which takes $[X] \in \C$ to the sum in $\F$ of the classes of $X\hbox{'s}$ weak orbits.  This map is well-defined by Property~3 above; it is a monoid homomorphism by Property~4; and it is surjective because $\phi([W]) = [W]$ if $W$ is a weak orbit.  It remains to show that $\phi$ is injective.

	Because $\F$ is free, there is a monoid homomorphism $\psi : \F \to \C$ which sends $[W]$ to $[W]$ for every weak orbit~$W$. It follows that $\psi \circ \phi$~is the identity on weak orbit classes.  But by Property~2 above, the weak orbit classes generate~$\C$, so in fact we have $\psi \circ \phi = \id_\C$. Since $\phi$ has a left inverse, it is injective, and hence a monoid isomorphism, as we wished to show.
\end{proof}

\subsubsection{Strong orbits}
	
	If $M = G$ is a group, then the structure of weak orbits is quite simple:  Indeed, if $x \ge y$, then $xg = y$ for some $g \in G$, hence $yg^{-1} = xgg^{-1} = x1 = x$, which means that $y \ge x$. Thus $\ge$ is symmetric and hence an equivalence relation.  It follows that a group action can take any element of a weak orbit to any other element.
		
	This is not true for monoids in general.  For example, let $M = \{0,1\}$ under multiplication, and let $M$ act on itself by right multiplication.  It is easy to see that $M$ is a weak orbit, but the action cannot take $0$ to~$1$.
		
	Define a relation $\sim$ on the elements of an $M$-set $X$ by $x \sim y$ if and only if $x \ge y$ and $y \ge x$. It is clear that $\sim$ is an equivalence relation, and we call its equivalence classes the \term{strong orbits} of $X$. If $X$ comprises a single equivalence class then we say that $X$ is a \term{strong orbit}.
		
	Note that strong orbits are contained in weak orbits.  As discussed above, weak orbits of a group action are strong orbits, so weak orbits and strong orbits coincide for a group action. This is usually not the case for monoid actions; for example, letting $M = \{0,1\}$ act on itself by right multiplication, we see that $\{1\}$ is a strong orbit but not invariant.
	
	It is easily verified that if:
		
	\[
	X \,=\, X_0 \,\supsetneq\, X_1 \,\supsetneq\, X_2 \,\supsetneq\, \cdots
	\,\supsetneq\, X_k \,\supsetneq\, X_{k+1} \,=\, \emptyset
	\]
	
\\ is an unrefinable chain of invariant subsets, then the strong orbits of~$X$ are $X_i \setminus X_{i+1}$ for $0 \le i \le k$.

	Finally, as with weak orbits, we may consider an auxiliary graph where $x \ge y$ is represented by a directed edge from~$x$ to~$y$. The strong orbits of $X$ are the strongly connected components of this graph.

\subsubsection{Properties of strong orbits}
	
	We now prove some elementary properties of strong orbits. 
	
\begin{Prop}\label{soprops}
The following hold:
\begin{enumerate}

\item The strong orbits of an $M$-set are subquotients.

\item Strong orbits of isomorphic $M$-sets can be put in a bijective, isomorphic correspondence.

\item A finite $M$-set contains an invariant strong orbit.

\item If $A$ is an invariant subset of an $M$-set $X$, then the strong orbits of $X$ are precisely those of $A$ and $X/A$, counting multiplicity.

\end{enumerate} 
\end{Prop}

\begin{proof}
	To show that a strong orbit $\Om$ is a subquotient, we must show that $\a mn \in \Om$ implies $\a m \in \Om$ for all $\a \in \Om$ and $m,n \in M$. Since $\a mn$ and $\a $ are in the same strong orbit, there exists $p \in M$ such that $\a = (\a mn)p = (\a m)(np)$, and hence $\a m \ge \a$. We always have $\a \ge \a m$, so $\a \sim \a m$, and hence $\a m \in \Om$, as desired.  This proves the first property.

	If $f$ is an isomorphism, then $x \ge y$ if and only if $f(x) \ge f(y)$. Since $\sim$~is defined in terms of $\ge$, it follows that the strong orbit of $x$ corresponds isomorphically via~$f$ to the strong orbit of $f(x)$. This proves the second property.
	
	Let $x$ be minimal with respect to $\ge$, and let $\Om$ be the strong orbit of~$x$. Then for any $\a \in \Om$ and $m \in M$ where $\a m$~is defined, we have $x \ge \a \ge \a m$. Minimality implies $x \sim \a m$, hence $\a m \in \Om$. Thus $\Om$ is invariant, which proves the third property.
	
	If $A$ is an invariant subset of $X$, then $a \in A$ and $a \ge x$ imply $x \in A$, hence $x \sim y$ implies that $x$ and~$y$ are both in~$A$ or both in $X\setminus A$. Since the actions on $A$ and $X/A$ are restrictions of the action on~$X$ to $A$ and $X \setminus A$, respectively, it follows that the strong orbits of $X$ are precisely the strong orbits of $A$ and~$X/A$. This settles the fourth property and concludes the proof.
\end{proof}

\subsection{The Burnside rings of a monoid}

	We are ready to explain the constructions which constitute the \emph{raison d'\^etre} of this paper.

\subsubsection{The weak Burnside ring}
	
	Earlier we showed that $(\C,+)$ is a commutative monoid, freely generated by the weak orbit classes.  Hence $(\C,+)$ embeds in $(\B_w,+)$, the free abelian group generated by the weak orbit classes.  Furthermore, since $(\C,+,\cdot)$ is a commutative semiring, it embeds in a commutative ring $(\B_w,+,\cdot)$, where $\cdot$ is defined on $\B_w$ in the natural way.  This is the Grothendieck-style completion alluded to in the Introduction.
	
	We call $\B_w$ the \term{weak Burnside ring} of~$M$.

\subsubsection{The strong Burnside relations}

	We wish to define a Burnside ring generated as an additive group by strong orbit classes, but the problem is that $M$-sets cannot be decomposed into strong orbits: strong orbits are subquotients, but not invariant in general.  To overcome this problem, we need to add relations to the weak Burnside ring which allow for such decompositions.
	
	Towards this end, let $I$ be the additive subgroup of $\B_w$ generated by elements of the form $[X]-[A]-[X/A]$, where $A$ is an invariant subset of $X$. We show that $I$ is an ideal.  Let $Y$ be an $M$-set, and observe that:
	
	\[
	([X]-[A]-[X/A]) \cdot [Y] \ = \
	[X \sprod Y] - [A \sprod Y] - [(X/A) \sprod Y] .
	\]
	
\\ Thus it suffices to show that $A \sprod Y$ is an invariant subset of $X \sprod Y$, and that:

	\cl $(*)$
	$(X/A) \sprod Y \, \cong\, (X \sprod Y) / (A \sprod Y)$. \\
	
	It is clear that $A \sprod Y$ is a subset of $X \sprod Y$. Moreover, if $(a,y)m$ is defined then it equals $(am,ym) \in A \sprod Y$, so $A \sprod Y$ is invariant, as desired.

	As for~$(*)$, notice that the underlying sets for the two $M$-sets both equal $(X \setminus A) \times Y$, so it suffices to show that the identity map is a morphism.  This in turn reduces to showing that:
	
	\[
	\hbox to 3.0em{$(x,y)m$\hfill} \qq\hbox{is defined in} \hskip1.5em (X/A) \sprod Y 
	\]
	
\\ if and only if:

	\[
	\hbox to 3.0em{$(x,y) \,\cdot\, m$\hfill} \qq\hbox{is defined in} \hskip1.5em (X \sprod Y) / (A \sprod Y).
	\]
	
\\ But this is easy to see, because both are defined precisely when $xm$ and $ym$ are defined and $xm \in X \setminus A$. Hence $I$ is an ideal of $\B_w$.

\subsubsection{The strong Burnside ring}

	We define $\B$ to be the quotient ring $\B_w/I$, and we call $\B$ the \term{strong Burnside ring} of~$M$. If we need to make $M$ explicit, we write $\B = \B(M)$. (The quotient $\B_w/I$ is the partial $M$-set analogue of what the authors of~\cite[Section 2]{weibel} call $K_0(\mathcal{C})$ in the $K$-theoretic, pointed $M$-set approach.)

	By an abuse of notation, we continue to write $[X]$ instead of $[X]+I$, and we disambiguate by indicating in which ring we are calculating.  In practice, this should not cause much grief, as we will not be discussing $\B_w$ beyond this section.
	
	We now prove that $\B$ is freely generated as an additive group by the strong orbit classes.  The following lemma may be compared with Property~2 of weak orbits, proved above in Section~2.3.2, Proposition~\ref{woprops}.

\begin{Lemma}\label{sogenerate}
If $X$ is a finite $M$-set, then in $\B$, $[X]$ equals the sum of the classes of $X$'s strong orbits.
\end{Lemma}

\begin{proof}
By Property~3 of strong orbits in Section~2.3.4, Proposition~\ref{soprops}, there exists an invariant strong orbit $\Om$ of~$X$. Thus we have $[X] = [\Om] + [X/\Om]$ in~$\B$. Property~4 of that same proposition implies that the strong orbits of~$X$ are precisely those of $\Om$ and $X/\Om$. By induction, we are done.
\end{proof}

\begin{Thm}
The strong orbit classes freely generate the additive group of $\B$.
\end{Thm}

\begin{proof}
Let $\F$ be a free abelian group on the strong orbit classes, and let $\phi : \B_w \to \F$ be the map which takes $[X] \in \B_w$ to the sum in $\F$ of the classes of $X$'s strong orbits.  This map is well-defined by Property~2 of strong orbits, and it is a group homomorphism by Property~4.

	Now suppose $A$ is an invariant subset of an $M$-set $X$. By Property~4 again, the strong orbit classes of $X$ are the strong orbit classes of $A$ and $X/A$. Hence $\phi([X]-[A]-[X/A]) = {\bf 0}$, so $I \setin \ker(\phi)$ and therefore $\phi$~factors through $\tilde\phi : \B \to \F$. For a strong orbit $\Om$, we have $\tilde\phi([\Om]) = [\Om]$, so $\tilde\phi$ is surjective.  We complete the proof by showing that $\tilde\phi$ is injective.

	Because $\F$ is free, there is a group homomorphism $\psi : \F \to \B$ which sends $[\Om]$ to $[\Om]$ for every strong orbit~$\Om$. It follows that $\psi \circ \tilde\phi$~is the identity on strong orbit classes.  But Lemma~\ref{sogenerate} shows that the strong orbit classes generate $\B$, hence we have $\psi \circ \tilde\phi = \id_\B$. Thus $\tilde\phi$ is injective, and a group isomorphism as desired. This completes the proof.
\end{proof}

\subsection{Final comments}

	Now that we have completed our constructions, let us take stock of what we have accomplished so far.
	
	We defined two Burnside rings, $\B_w$ and $\B$, whose additive groups are freely generated by the weak and strong orbit classes, respectively.  But with group actions, weak and strong orbits coincide, so $\B_w \cong \B$, and hence there is only one Burnside ring of a group.

	We will see that $\B$, rather than $\B_w$, is the appropriate correlate of the Burnside ring of a group. For instance, it is well-known that every transitive $G$-set is isomorphic to a system of cosets $G/K$, where $K$ is a subgroup of $G$, and that $G/K \cong G/L$ precisely when $K$ and $L$ are conjugate.  Therefore the structure of the Burnside ring of a group~$G$ is intimately connected with the subgroup structure of~$G$. We will prove similar results for strong orbits in the coming sections.  No such results are possible for weak orbits.
	
	In fact, weak orbits are quite unruly, as a simple example will show: Consider the monoid $M = \{0,1\}$ under multiplication.  For every natural number~$n$, let $X_n = \{0,1,\dots,n\}$. Then $M$ acts on $X_n$ by right multiplication.  Observe that each $X_n$ is a weak orbit, but distinct $X_n$ are not isomorphic, since they have different numbers of elements.  This means that even for a two-element monoid, the weak Burnside ring may have infinite rank!
	
	By contrast, the only strong orbits of~$M$ are $\{0\}$, which is total and invariant, and $\{1\}$, which is only a subquotient, since $1 \rac 0 = 0 \notin \{1\}$. It is easy to see that $[\{0\}] = {\bf 1}$, so letting $\omega = [\{1\}]$, we have:
	
	\[
	[X_n] \,=\, {\bf 1} + n\omega \qq\hbox{in $\B$,} 
	\]
	
\\ a decomposition which neatly captures the structure of this family of $M$-sets. In fact, one can straightforwardly check that $\omega \cdot \omega = \omega$, and hence that $\B \cong \Z \times \Z$ as rings via the map $m{\bf 1}+n\omega \,\mapsto\, (m,m+n)$.

	Therefore, going forward we focus on $\B$ exclusively, which shall henceforth be known simply as \term{the Burnside ring of a monoid}.

\section{Strong orbits and the structure of finite monoids}

	Every transitive right $G$-set is isomorphic to a system of right cosets $G/K$ under the action of right multiplication, where $K$ is a subgroup of $G$; and the automorphism group of $G/K$ is isomorphic to $\mathcal{N}_G(K)/K$, where $\mathcal{N}_G(K)$ is the normalizer of~$K$ in~$G$. In this section, we give the corresponding theory for finite monoids, due essentially to Green~\cite{green} and Rhodes~\cite{rhodes}.
	
	These results do not hold in general for infinite monoids, so we restrict our attention in the sequel to finite monoids.

\subsection{The action of right multiplication}

	Any monoid acts on itself by \term{right multiplication}: in symbols, $x \rac m = xm$.
		
	Recall that we defined the relation $\ge$ on the elements of an $M$-set by stipulating that $x \ge y$ holds if and only if there exists some $m \in M$ such that $x \rac m = y$.  We defined an equivalence relation $\sim$ where $x \sim y$ if and only if $x \ge y$ and $y \ge x$, and the equivalence classes of~$\sim$ are the strong orbits.
	
	In the context of right multiplication, we write $\ger$ for $\ge$ and $\R$ for $\sim$. Thus $x \ger y$ holds if and only if there exists $m \in M$ such that $xm = y$, and $x \R y$ holds if and only if $x \ger y$ and $y \ger x$. The strong orbits are the \term{$\R$-classes} $R_x = \{ y \in M : x \R y \}$.
	
	Recall that a congruence on an $M$-set $X$ is an equivalence relation $\equiv$ such that $x \equiv y$ implies $x \rac m \equiv y \rac m$ for all $x,y \in X$ and $m \in M$. In the context of right multiplication, this means $x \equiv y$ implies $xm \equiv ym$ for all $x,y,m \in M$. We call such a relation a \term{right congruence} on $M$.
	
	Every strong orbit of $M$ is isomorphic to a quotient of an $\R$-class by a right congruence, a result which will be proved below in Proposition~\ref{apex}.

\subsection{Finite monoid theory}

	We introduce some standard concepts from finite monoid theory.  We will also borrow a few facts about monoids in the sequel.  Proofs and a more thorough discussion can be found in \cite{qtheory}.

\subsubsection{Green's relations}

	To the relations $\ger$ and $\R$ we add \term{Green's preorders}:

	\medbreak
	\ul \hbox to 4.0em{$x \gel y$\hfil} $\Longleftrightarrow \q$ \hbox to 4.0em{$mx = y$\hfil} for some $m \in M$\\
	\ul \hbox to 4.0em{$x \gej y$\hfil} $\Longleftrightarrow \q$ \hbox to 4.0em{$mxn = y$\hfil} for some $m,n \in M$\\
	\ul \hbox to 4.0em{$x \geh y$\hfil} $\Longleftrightarrow \q$ $x \gel y$ \ {\rm and} \ $x \ger y$\\
	\medbreak

\\ and their corresponding equivalence relations $\L$, $\J$, and $\H$, known as \term{Green's relations}. Just as we write $R_x$ for the $\R$-class of~$x$, we write $L_x$, $J_x$, and $H_x$ for the equivalence classes of~$x$ under the appropriate relation.  Note that the containments:

	\[
	\H \ \setin \ \L,\R \ \setin \ \J
	\]
	
\\ are immediate from the definitions.  A group comprises a single $\H$-class.

	If $x \ger y$, we say that $x$ is \term{$\R$-above} $y$, and similarly for other relations.  If $x$ is $\J$-above~$y$, then any element of $J_x$ is $\J$-above any element of~$J_y$.  In this way, we may speak of relations holding between classes as well as elements.

	Finite monoids satisfy a property called \term{stability} (first introduced as \term{weak stability} in \cite{ocarroll}): if $x \J y$, then $x \ger y$ implies $x \R y$ and $x \gel y$ implies $x \L y$. \cite[Theorem A.2.4]{qtheory} We often use stability in the following formulation: Let $J$ be the $\J$-class of $x$, and multiply $x$ on the right by~$m$. If $xm \in J$ as well, then $x \R xm$ by stability. Similarly for multiplication on the left.

\subsubsection{Idempotents}

	We say that $e \in M$ is an \term{idempotent} if $ee = e$. It is easy to check that if $e$ is an idempotent of $M$, then $e$~is a right identity for~$L_e$, a left identity for~$R_e$, and a two-sided identity for $H_e$. \cite[Lemma~4]{cliffordmiller} We will use these properties frequently in the sequel.

	A \term{subgroup} is a subset of $M$ which forms a group under the operation of~$M$, not necessarily with $1$ as the identity.  The $\H$-class of an idempotent~$e$ is a subgroup of $M$ with identity~$e$. \cite[Theorem 7]{green} In fact, it is the group of units of the monoid $eMe$, and contains every subgroup of~$M$ with identity $e$. For this reason, $H_e$ is called a \term{maximal subgroup}.

\subsection{The apex of a strong orbit}

	Let $\Omega$ be a strong orbit of~$M$. We say that $m \in M$ \term{annihilates} $\a \in \Omega$ if $\a m \notin \Omega$; equivalently, that $\a \cdot m$ is undefined.  Similarly, we say that $m$ annihilates $\Omega$ if $m$ annihilates every element of $\Omega$.
	
	The following proposition is fundamental.  We include a proof to collect in one place the particular properties we use here; another approach can be found in Margolis and Steinberg~\cite[Lemma 2.1]{margolis}.
		
\begin{Prop}\label{apex}
Let $\Omega$ be a strong orbit.  The following hold: 
\begin{enumerate}

\item There exists a unique $\J$-class $J$, called the \term{apex} of $\Omega$, such that $m \in M$ does not annihilate~$\Omega$ if and only if $m$ is $\J$-above $J$.

\item For all $\a \in \Omega$, there exists an idempotent $e \in J$ such that $\a e = \a$ and $\a r$ is defined for all $r \in R_e$. Furthermore, $\Omega = \a R_e$ and $\Omega$ is isomorphic to a quotient $R_e/{\equiv}$.

\item Isomorphic strong orbits have the same apex.

\item For any idempotent $e$, the apex of $R_e$, or a quotient $R_e/{\equiv}$, is $J_e$.

\end{enumerate}
\end{Prop}

\begin{proof}

	Fix $\a \in \Omega$. Let $e$ be an $\R$-minimal idempotent which fixes~$\a$, and let $J$~be the $\J$-class of~$e$. (This is always possible since $1$ fixes $\a$ and $M$ is finite.) If $m$~does not annihilate~$\Omega$, then $\b m \in \Omega$ for some $\b \in \Omega$, and since $\Omega$ is a strong orbit, there exist $x,y \in M$ such that $\a x = \b$ and $\b my = \a$, whence $\a(exmy) = \a$. Since $M$ is finite, there is a positive integer $k$ such that $(exmy)^k$ is an idempotent \cite[Section 1.2]{rees}, hence $(exmy)^k$ is an idempotent which fixes~$\a$. By $\R$-minimality of $e$ we have $(exmy)^k \ger e$, hence $m \gej e$. Conversely, if $m \gej e$ then $\a xmy = \a \in \Omega$ for some $x,y \in M$. Hence $\a x$ and $(\a x)m$ are in~$\Omega$, so $m$ does not annihilate~$\Omega$. It follows that $J$~has the required property, and uniqueness is immediate. This settles the first property.
		
	If $e \ger m$ and $m$ does not annihilate~$\a$, then $m \gej e$ and hence $e \R m$ by stability. This shows that $e$ is an $\R$-minimal element which does not annihilate~$\a$. Now let $r \in R_e$ and $x \in M$. If $rx \in R_e$ then $rxy = e$ for some $y \in M$, hence $\a rxy = \a \in \Omega$, so $\a rx \in \Omega$. Conversely, if $\a rx = \a erx\in \Omega$ then $rx \in R_e$ by minimality. Thus we have shown that $rx \in R_e$ if and only if $\a rx \in \Omega$, hence $r \mapsto \a r$ is a well-defined morphism from $R_e$ to~$\Omega$. Additionally, if $\a x = \a ex \in \Omega$, then $ex \in R_e$, so this morphism is surjective.  Hence $\Omega \cong R_e/{\equiv}$, where $r \equiv s$ if and only if $\a r = \a s$. This settles the second property.
	
	If $f$ is an isomorphism of strong orbits, then $m$ annihilates $\a$ if and only if $m$ annihilates $f(\a)$, whence the third property follows.

	Finally, let $e$ be an idempotent and let $J_x$ be the apex of $R_e/{\equiv}$, where $x$ does not annihilate $\overline e$. Then $\overline e \cdot x = \overline{ex}$, so $ex \in R_e$ and hence $x \gej e$. But $\overline e \cdot e = \overline{ee} = \overline e$, so $e$ does not annihilate $R_e/{\equiv}$, and hence $e \gej x$. Thus $x \J e$ and $J_x = J_e$ as claimed. Since $R_e \cong R_e/{=}$ and isomorphic strong orbits have the same apex, it follows that $J_e$ is the apex of $R_e$ as well.
\end{proof}

	We will use the characterizations of strong orbits in Property 2 quite frequently in the sequel.
	
	A finite monoid has finitely many $\R$-classes, and since these classes are also finite, there are finitely many quotients of these.  Hence Property~2 implies that the additive group of the Burnside ring of a finite monoid is free of finite rank.

	The same need not be true for an infinite monoid: Consider the free monoid on a single letter~$a$. Each word forms a one-element $\R$-class which is annihilated by every word except the empty word. But for every positive integer~$n$, there is a cyclic strong orbit $\{x_0,\dots,x_{n-1}\}$ with $x_i \rac a = x_{(i+1) \, {\rm mod} \, n}$ which has no annihilators. Hence for infinite monoids, strong orbits are not necessarily given by quotients of $\R$-classes, and the additive group of the Burnside ring may be free of infinite rank.

\subsection{Automorphisms of a strong orbit}

	We conclude this section by computing the group of automorphisms of a strong orbit.  Recall that in the group setting, the automorphism group of a system of right cosets $G/K$ is isomorphic to $\mathcal{N}_G(K)/K$. We may prove this by first observing that $\mathcal{N}_G(K)$ acts on the left of $G/K$ by multiplication, with kernel $K$. We then check that every automorphism of $G/K$ is given by left multiplication by an element of $\mathcal{N}_G(K)$. A similar approach will work for monoids in general.
	
\subsubsection{Left multiplication and automorphisms of $R_e$}
		
	Let $\l_m : M \to M$ be the left multiplication map, defined by $\l_m(x) = mx$, and let $e$~be an idempotent of~$M$. The group $H_e$ acts freely on $R_e$ by left multiplication; in symbols, $h \lac r = \l_h(r) = hr$ for all $h \in H_e$ and $r \in R_e$. The action is by automorphisms, meaning that $\l_h(r \cdot m) = \l_h(r) \cdot m$ for all $m \in M$.
	
	Thus each element $h \in H_e$ yields an automorphism $\l_h$ of $R_e$, and in fact it is known that $\Aut(R_e) \cong H_e$. In the next section, we generalize this result by determining more generally the structure of $\Aut(R_e/{\equiv})$.
	
\subsubsection{Automorphisms of $R_e/{\equiv}$}
	
	Recall that by Proposition~\ref{apex}, every strong orbit is isomorphic to a quotient $R_e/{\equiv}$, where $e$~is an idempotent and $\equiv$ is a right congruence.
	
	Consider a strong orbit $R_e/{\equiv}$, and let $\fL$~consist of all elements of $H_e$ which preserve $\equiv$ on the left; in symbols, $h \in \fL$ if and only if $r \equiv s$ implies $hr \equiv hs$ for all $r,s \in R_e$. Since $e$~is a left identity for~$R_e$, we have $e \in \fL$. Also, it is clear by associativity that $\fL$ is closed under multiplication.  Since $H_e$ is finite, we conclude that $\fL$ is a subgroup of~$H_e$. It follows that $\fL$ acts on the left of~$R_e/{\equiv}$ by automorphisms, where we define $h \lac \overline r = \overline{hr}$ for all $h \in \fL$ and $r \in R_e$.
	
	If $h \lac \overline r = \overline r$ for some $r \in R_e$, then $h \lac (\overline r \cdot m) = (h \lac \overline r) \cdot m = \overline r \cdot m$ for all $m \in M$, hence $h$~acts trivially on $\overline r M$, which contains $R_e/{\equiv}$.  Thus, letting $\fK$~be the kernel of the action, we see that $h \in \fK$ if and only if $h$~acts trivially on one element of $R_e/{\equiv}$.  We conclude that $\fL/\fK$ acts freely on~$R_e/{\equiv}$, and so there exists an injective group homomorphism $\phi : \fL/\fK \to \Aut(R_e/{\equiv})$ defined by $\phi(\fK \ell) = \Lambda_\ell$, where $\Lambda_\ell(\overline r) = \overline{\ell r}$.
	
	Finally, we show that $\phi$ is surjective.  Let $f$ be an arbitrary automorphism of $R_e/{\equiv}$, and suppose that $f(\overline e) = \overline r$ for $r \in R_e$. We show that $re \in \fL$ and $f = \Lambda_{re} = \phi(\fK re)$.
	
	Recall the principle mentioned in Section~1.2: if we begin a calculation with a defined expression, then all subsequent expressions are defined as well. Thus observe that:
	
	\[
	\overline r \,=\, f(\overline e) \,=\, f(\overline {ee}) \,=\, f(\overline e \cdot e) \,=\, f(\overline e) \cdot e \,=\, \overline r \cdot e \,=\, \overline{re} 
	\]
	
\\ so $re \in R_e$ and $\overline r = \overline{re}$. But $re \in R_e$ implies $re \in J_e$, whence $re \in L_e$ by stability, and thus $re \in H_e$.

	Now let $x \in R_e$ be arbitrary.  Then we have:
	
	\[
	f(\overline x) \,=\, f(\overline{ex}) \,=\, f(\overline e \cdot x) \,=\, f(\overline e) \cdot x \,=\, \overline{re} \cdot x \,=\, \overline{rex} 
	\]
	
\\ so $rex \in R_e$ and $f(\overline x) = \overline{rex}$. But then if $y \in R_e$ with $x \equiv y$, we have:

	\[
	\overline{rex} \,=\, f(\overline x) \,=\, f(\overline y) \,=\, \overline{rey} 
	\]
	
\\ so $rex \equiv rey$. Thus $re \in \fL$ and $f(\overline x) = \overline{rex} = \Lambda_{re}(\overline x)$ as claimed.
	
	We have therefore proved the following:

\begin{Thm}\label{autstructure}
$\Aut(R_e/{\equiv})$ acts freely and is isomorphic to $\fL/\fK$. \qed
\end{Thm}

	We conclude with a few observations.
	
	First, if we take $=$ for $\equiv$, then $R_e/{\equiv}$ is~isomorphic to $R_e$, $\fL$ is $H_e$, and $\fK$ is trivial, so our theorem yields the well-known result $\Aut(R_e) \cong H_e$ mentioned above.
	
	Second, as we saw above, we have $h \in \fK$ if and only if $\overline h = h \lac \overline e = \overline e$, hence we have $\fK = \{ h \in \fL : h \equiv e \}$. But if $h \in H_e$ and $h \equiv e$, then $hr \equiv er = r$, so $r \equiv s$ implies $hr \equiv r \equiv s \equiv hs$, and hence $h \in \fL$, so in fact we have $\fK = \{ h \in H_e : h \equiv e \}$.

	In particular, if $M = G$ is a group, then it is a standard result that right congruence classes are right cosets of the class of the identity, so $R_e/{\equiv}$ is $G/\fK$. It is easily verified that $\fL$~is $\mathcal{N}_G(\fK)$, hence our theorem yields the familiar result $\Aut(G/\fK) \cong \mathcal{N}_G(\fK)/\fK$. In light of this, the reader might reasonably wonder whether $\fL = \mathcal{N}_{H_e}(\fK)$ holds generally.  However, Benjamin Steinberg constructed a counterexample, which can be found as an Appendix.

\section{Distinguishability}

	In this section, we define \term{regular $\J$-classes} and the notion of \term{distinguishability}, the latter of which has been explored recently by Margolis and Steinberg in~\cite[Proposition 2.9]{margolis}. We show that a monoid is distinguishable if and only if every strong orbit is isomorphic to a system of cosets of a maximal subgroup~$H_e$, and that a monoid is distinguishable if and only if its Burnside ring is isomorphic to the direct product of the Burnside rings of maximal subgroups. We conclude with several examples. 
	
	We continue to restrict our attention to finite monoids.
	
\subsection{Definitions}

	We say that a $\J$-class is \term{regular} if it contains an idempotent.  Thus, for example, the apex of a strong orbit is regular.
	
	Let $J$ be a regular $\J$-class. We say that $\L$-classes $L$ and~$L'$ in~$J$ are \term{distinguishable} if there exists an $\R$-class~$R$ in~$J$ such that precisely one of $L \cap R$ and $L' \cap R$ contains an idempotent.  In this case, we say that $R$ \term{distinguishes} $L$ and~$L'$. If no such $\R$-class exists then we say that $L$ and~$L'$ are \term{indistinguishable}. If $L = L'$ then $L$ and~$L'$ are trivially indistinguishable.
	
	We say that $J$ itself is \term{distinguishable} if every pair of distinct $\L$-classes in~$J$ is distinguishable.  Finally, we say that $M$~is \term{distinguishable} if every regular $\J$-class in~$M$ is distinguishable. A group comprises a single $\L$-class, and is hence trivially distinguishable.

	Indistinguishability can be helpfully rephrased using the Clifford-Miller theorem~\cite[Theorem~3]{cliffordmiller}:

\begin{Thm}[Clifford-Miller]\label{cmthm}
Let $J$ be a $\J$-class and let $x,y \in J$. Then $xy \in J$ if and only if $L_x \cap R_y$ contains an idempotent. \qed
\end{Thm}

\\ Hence if $x$ and $x'$ are elements of a regular $\J$-class $J$, then $L_x$ and $L_{x'}$ are indistinguishable if and only if:

	\cl $(*)$
	$xy \in J \ \Longleftrightarrow\ x'y \in J$ \qq\hbox{for all $y \in J$.}\\
	
\\ In fact, in this case $(*)$ holds for all $m \in M$: For example, if $xm \in J$ then by Proposition~\ref{apex} there is an idempotent $y \in J$ such that $xmy = xm \in J$. It follows that $my \in J$, so $x'my \in J$ by~$(*)$, and hence $x'm \in J$, as claimed.

\subsection{Congruences on $R_e$}

	Recall from Section~3 that every strong orbit is isomorphic to a quotient $R_e/{\equiv}$, where $R_e$~is the $\R$-class of an idempotent, and $\equiv$~is a right congruence. We now discuss a particular congruence relation which is determined by a subgroup of $H_e$. We show that if the surrounding $\J$-class is distinguishable, then all right congruences on $R_e$ are of this type.

\subsubsection{Quotients and right congruences determined by maximal subgroups}

	Let $e$~be an idempotent, and let $R_e$ be the $\R$-class of $e$. We saw in Section~3 that $H_e$ acts on the left of $R_e$ by automorphisms, hence so does any subgroup $K \le H_e$. Because this left action is by automorphisms, it follows that $M$ acts on the right of the set of orbits $\{Kr : r \in R_e\}$ by letting $Kr \rac m = Krm$ if $rm \in R_e$ and letting $Kr \rac m$ be undefined otherwise.

	Now define $\equiv_K$ on $R_e$ by $r \equiv_K s$ if and only if $Kr = Ks$. It follows from the above that $\equiv_K$ is a right congruence.  We call $\equiv_K$ the right congruence on~$R_e$ \term{determined by $K$}. More generally, we say such congruences are \term{determined by maximal subgroups}. We write $R_e/K$ for $R_e/{\equiv}_K$, and it is clear that $R_e/K = \{ Kr : r \in R_e \}$.
	
	 If a strong orbit is isomorphic to $R_e/K$, we say that the strong orbit is \term{determined by $K$}. More generally, we say that a strong orbit of this form is \term{determined by maximal subgroups}.

\subsubsection{Right congruences contained in $\H$}
	
	The right congruence $\equiv_K$ is contained in $\H$, meaning that $r \equiv_K s$ implies $r \H s$; or, put another way, that the congruence classes of $\equiv_K$ are contained in $\H$-classes. Indeed, if $r \equiv_K s$, then $Kr = Ks$, so $r = ks$ and $s = k'r$ for some $k,k' \in K$. This shows that we have $r \L s$. But we also have $r \R s$, since $r,s \in R_e$, hence we have $r \H s$, as claimed.

	Conversely, we show that any right congruence on $R_e$ contained in $\H$ is determined by a subgroup of $H_e$, namely $\overline e$, the congruence class of~$e$.
	
	Let $\equiv$ be a right congruence on $R_e$ contained in $\H$. Then $\overline e$ is a subset of~$H_e$ which contains~$e$. Moreover, for any $r,s \in \overline e$, we have:
	
	\[
	e \,\equiv\, s \,=\, es \,=\, e \cdot s \,\equiv\, r \cdot s \,=\, rs 
	\]
	
\\ which shows that $\overline e$ is closed under multiplication.  Similarly, we have:

	\[
	r^{-1} \,=\, er^{-1} \,=\, e \cdot r^{-1} \,\equiv\, r \cdot r^{-1} \,=\, rr^{-1} \,=\, e 
	\]
	
\\ which shows that $\overline e$ is closed under taking inverses, so $\overline e$ is a subgroup of~$H_e$.

	It remains to show that $r \equiv s$ if and only $r \equiv_{\,\overline e} s$. First suppose $r \equiv_{\,\overline e} s$. Then there exists $k \in \overline e$ such that $kr = s$. Thus we have:
	
	\[
	r \,=\, er \,=\, e \cdot r \,\equiv\, k \cdot r \,=\, kr \,=\, s 
	\]
	
\\ as claimed.

	Conversely, suppose $r \equiv s$. Since $\equiv$ is contained in~$\H$, we have $r \H s$, so there exists $m \in M$ such that $r = ms$. Because $s \in R_e$, we may also write $e = sn$ for some $n \in M$. Thus we have:
	
	\[
	e \,=\, sn \,=\, s \cdot n \,\equiv\, r \cdot n \,=\, rn \,=\, msn \,=\, me 
	\]
	
\\ and so $me \in \overline e$. Since $r = ms = m(es) = (me)s$, we have $r \equiv_{\,\overline e} s$ as claimed.

	Thus we have established the following:
	
\begin{Thm}\label{congchar}
Right congruences contained in $\H$ are precisely the right congruences determined by maximal subgroups.  Specifically, a right congruence on~$R_e$ contained in $\H$ is determined by $\overline e \le H_e$. \qed
\end{Thm}

\subsubsection{Right congruences and distinguishability}

	Let $R_e$ be the $\R$-class of an idempotent, and let $\equiv$ be a right congruence on~$R_e$.

\begin{Prop}\label{congindis}
Let $r$ and $s$ be elements of $R_e$. If $r \equiv s$, then $L_r$ and $L_s$ are indistinguishable.
\end{Prop}

\begin{proof}
If $r \equiv s$, then $r \cdot y \equiv s \cdot y$ for all $y \in J_e$, hence $ry \in J_e$ if and only if $sy \in J_e$. Thus by Clifford-Miller, $L_r$ and $L_s$ are indistinguishable.
\end{proof}

\begin{Cor}\label{distcong}
If $J_e$ is distinguishable, then $\equiv$ is contained in $\H$, and hence $\equiv$ and $R_e/{\equiv}$ are determined by maximal subgroups.
\end{Cor}

\begin{proof}
Suppose $r \equiv s$. Then $L_r$ and $L_s$ are indistinguishable by the preceding proposition. Since $J_e$ is distinguishable, we must have $r \L s$. Since $r,s \in R_e$, we already have $r \R s$, so $r \H s$. The final claim follows from Theorem~\ref{congchar}.
\end{proof}

\begin{Cor}\label{burnsideg}
If $M$ is distinguishable, then the additive subgroup of $\B$ is generated by the classes of strong orbits which are determined by maximal subgroups.
\end{Cor}
\begin{proof}
Every strong orbit is isomorphic to $R_e/{\equiv}$, where $e$ is an idempotent, and hence $J_e$~is regular.  If $M$~is distinguishable, then $J_e$ is distinguishable, and hence $R_e/{\equiv}$ is determined by maximal subgroups by the preceding corollary.  The conclusion follows because the additive subgroup of $\B$ is generated by the strong orbit classes by Lemma~\ref{sogenerate}. \end{proof}

\subsubsection{The maximal right congruence on $R_e$}

	Once again, let $R_e$ be the $\R$-class of an idempotent, and define $\emax$ on $R_e$ by $r \emax s$ if and only if $L_r$ and~$L_s$ are indistinguishable.  We show that $\emax$ is a right congruence on $R_e$. 
	
	For $r,s \in R_e$, recall that $L_r$ and $L_s$ are indistinguishable if and only if:
	
	\cl $(*)$
	$rm \in J_e \ \Longleftrightarrow\ sm \in J_e$ \qq\hbox{for all $m \in M$.}\\
	
\\ Since $\Leftrightarrow$ is an equivalence relation, so is $\emax$. Now let $r \emax s$ hold for $r,s \in R_e$ and let $m \in M$ be arbitrary.  If $rm \notin J_e$ and $sm \notin J_e$, then $r \cdot m$ and $s \cdot m$ are both undefined, so $r \cdot m \emax s \cdot m$ holds.  Otherwise, by $(*)$ we have $rm$ and~$sm$ in~$J_e$, hence in~$R_e$ by stability. For any $n \in M$, apply $(*)$ again with~$mn$ to conclude that $rmn \in J_e$ if and only if $smn \in J_e$. Since $n$~is arbitrary, $L_{rm}$ and $L_{sm}$ are indistinguishable by~$(*)$, and hence $r \cdot m \emax s \cdot m$.
	
	Thus $\emax$ is a right congruence.  Since elements from distinguishable $\L$-classes cannot be congruent in any right congruence by Proposition~\ref{congindis}, it follows that $\emax$ is the maximal right congruence on~$R_e$.

	When $J_e$ is distinguishable, $\emax$ is $\equiv_{H_e}$.  Indeed, if $J_e$~is distinguishable, then $\emax$ is $\equiv_{\,\overline e}$ by Corollary~\ref{distcong}.  But $h \in H_e$ implies $L_h = L_e$, and hence $h \emax e$. Thus $h \in \overline e$, and so $\overline e = H_e$, as claimed.
	
	On the other hand, if $J_e$ is not distinguishable, then $\emax$ is not determined by maximal subgroups at all.  Indeed, for all $r,s \in R_e$, $r \H s$ implies $L_r = L_s$, so $r \emax s$ holds trivially.  Hence the $\H$-classes of $R_e$ are contained in the classes of~$\emax$, and so:

	\cl $(\dag)$
	$|R_e/{\emax}| \ \le \ \hbox{the number of $\H$-classes in $R_e$}$.\\

\\ Now let $f$ be an idempotent and let $K \le H_f$.  Then the classes of~$\equiv_K$ are contained in the $\H$-classes of~$R_f$, and so:

	\cl $(\ddag)$
	$\hbox{the number of $\H$-classes in $R_f$} \ \le \ |R_f/K|$.\\

\\ If $R_e/{\emax} \cong R_f/K$, then $e \J f$ by considering the apex, and hence by a standard result of monoid theory, $R_e$ and~$R_f$ have contain the same number of $\H$-classes.  \cite[Corollary~A.3.2]{qtheory} But we also have $|R_e/{\emax}| = |R_f/K|$. This forces equality in~$(\dag)$, so the classes of $\emax$ are just the $\H$-classes of $R_e$.

	Finally, consider two $\L$-classes in $J_e$. By another standard result of monoid theory, these $\L$-classes intersect~$R_e$, so we can write them as $L_r$ and~$L_s$ for some $r,s \in R_e$. \cite[Corollary~A.2.5]{qtheory} If $L_r$ and $L_s$ are indistinguishable, then $r \emax s$, hence $r \H s$ by the above, and so $L_r = L_s$. Hence $J_e$ is distinguishable.
	
	Therefore if $J_e$ is not distinguishable, then $R_e/{\emax}$ is not isomorphic to any quotient of the form $R_f/K$, which means that $R_e/{\emax}$ is not determined by maximal subgroups, as claimed.

\subsection{The Burnside ring of a distinguishable monoid}

	In this section, we show that the Burnside ring of a distinguishable monoid is isomorphic to the direct product of the Burnside rings of its maximal subgroups.

\subsubsection{Definition of the isomorphism}

	In order to define our proposed isomorphism, we need a lemma.

\begin{Lemma}
Let $X$ be an $M$-set and let $e \in M$ be an idempotent.  Then:

	\[
	Xe \, = \, \{ \, xe \, : \, \hbox{$x \in X$ and $\,xe$ is defined} \hskip2.5pt \}
	\]
	
\\ is a right $H_e$-set, under the restriction of the original action to $H_e$.
\end{Lemma}
\begin{proof}
If $xe \in Xe$, then $xe$ is defined. But for all $h \in H_e$:

	\[
	xe \,=\, xehh^{-1}
	\]
	
\\ so $xeh$ is defined and $xeh \in X$.  Thus $xeh = xehe \in Xe$, and hence $Xe$ is closed under the action of~$H_e$. Associativity follows from associativity of the original action, and $(xe)e = x(ee) = xe$ shows that the action of the identity is correct.  Thus $Xe$ is an $H_e$-set, as claimed.
\end{proof}

	Now let $\idrep$ contain one idempotent from each regular $\J$-class in $M$, and let $e$ range over the elements of $\idrep$. We define:
	
	\[
	\phi : \B(M) \longrightarrow \prod_{e \in \idrep} \B(H_e)
	\]
	
\\ by:

	\[
	\phi([X]) = (e: [Xe] ) 
	\]
	
\\ where the notation $(e: P(e))$ stands for the tuple with $P(e)$ in the coordinate corresponding to~$e$. It is clear that $X \cong Y$ implies $Xe \cong Ye$, so the definition of $\phi$ is independent of representatives.

	Let $\G(M)$ by the additive subgroup of $\B(M)$ generated by the classes of strong orbits which are determined by maximal subgroups.  Recall that by Corollary~\ref{burnsideg}, if $M$~is distinguishable then $\B(M)=\G(M)$. We will show that $\phi$ is a surjective ring homomorphism that restricts to an group isomorphism on~$\G(M)$. (In general, $\G(M)$ is not a subring of~$\B(M)$.  For example, let $M = \{a,b,e,f,1\}$, all idempotents, where $ab=b$, $ba = a$, $ef = f$, $fe = e$, and where $e$ and~$f$ are identity elements for $\{a,b\}$.  Note $R_a = \{a,b\} \in \G(M)$. It can be checked that $R_a \sprod R_a$ contains a one-element strong orbit $\{(a,b)\}$ whose apex is $J_e = R_e = \{e,f\}$.  But $\{(a,b)\} \not\cong R_e/K$ for $K \le H_e$, because $H_e = \{e\}$ is trivial and so $R_e/K$ has two elements.)
	
	We first show that $\phi$ is a homomorphism.  Additivity follows from:
		
	\[
	(X \ssum Y)e \,=\, Xe \ssum Ye .
	\]
	
\\ To see this equality, recall that $X \ssum Y$ is the disjoint union of $X$ and $Y$. Thus $(X \ssum Y)e$ is the disjoint union of $Xe$ and $Ye$, which is the righthand side.  Similarly, multiplicativity follows from:

	\[
	(X \sprod Y)e \,=\, Xe \sprod Ye .
	\]
	
\\ To see this equality, recall that the elements of $X \sprod Y$ are the pairs $(x,y)$ with $x \in X$ and $y \in Y$, hence the elements of $(X \sprod Y)e$ are the pairs $(xe,ye)$ where $xe$ and $ye$ are defined.  These are exactly the elements on the righthand side.  Finally, we check that $\phi$ preserves the identity.  Recall that ${\bf 1} = [\{\s\}]$, where $\s m = \s$ for all $m \in M$. Thus:

	\[
	\phi({\bf 1}) \,=\, \phi([\{\s\}]) \,=\,
	(e: [\{\s\}e]\,) \,=\, (e: [\{\s\}]\,) \,=\, (e: {\bf 1}) 
	\]
	
\\ which is the identity of $\prod_e \B(H_e)$.
	
\subsubsection{Bases}

	In this section, we fix bases for $\G(M)$ and the additive group of $\prod_e \B(H_e)$ and show that they have the same size.
	
	Suppose $u$ and $v$ are idempotents of $M$ with $u \J v$. It is a standard result of monoid theory that $R_u \cong R_v$, and hence the quotients of $R_u$ and~$R_v$ are in a bijective, isomorphic correspondence. \cite[Lemma,~p.~165]{green} It follows that every strong orbit of~$M$ is isomorphic to a quotient of some $R_e$, $e \in \idrep$.

	We now show that for strong orbits determined by subgroups of~$H_e$, it suffices to take one subgroup from each conjugacy class.

\begin{Prop}\label{morphcond}
Let $K$ and $L$ be subgroups of $H_e$. Then there is a morphism from $R_e/K$ to $R_e/L$ if and only if $K$ is contained in a conjugate of $L$ by an element of $H_e$. More precisely, if $f : R_e/K \to R_e/L$ is a morphism satisfying $f(K) = Lr$, then $re \in H_e$ and $K \setin (re)^{-1}L(re)$; conversely, if $K \setin h^{-1}Lh$ for some $h \in H_e$ then there is a morphism $f$ satisfying $f(K) = Lh$.
\end{Prop}
\begin{proof}
If $f$ is a morphism satisfying $f(K) = Lr$, then for all $k \in K$ we have:

	\[
	Lr \,=\, f(K) \,=\, f(Kk) \,=\, f(K \rac k) \,=\, f(K) \rac k \,=\, Lr \rac k \,=\, Lrk 
	\]
	
\\ so $rk \in R_e$ and $Lr = Lrk$ for all $k \in K$. Taking $k = e$, we see that $re \in R_e$. Then $re \in J_e$, so $re \in L_e$ by stability, and therefore $re \in H_e$. Moreover, for any $k \in K$ we have:

	\[
	(re)k \,=\, r(ek) \,=\, rk \,=\, erk \,\in\, Lrk \,=\, Lr \,=\, Lre 
	\]
	
\\ or $k \in (re)^{-1}L(re)$. Thus $K \setin (re)^{-1}L(re)$, as claimed.

	Conversely, suppose $K \setin h^{-1}Lh$. We wish to define a map $f : R_e/K \to R_e/L$ by $f(Kr) = Lhr$. Recall that $H_e$ acts on the left of $R_e$, hence $hr \in R_e$ and so $Lhr \in R_e/L$. To show independence of representatives, we require $f(Kr) = f(Kkr)$, or $Lhr = Lhkr$, for all $k \in K$. Indeed, by hypothesis, for all $k \in K$ there exists $\ell \in L$ such that $k = h^{-1}\ell h$, and hence:
	
	\[
	Lhkr \,=\, Lh(h^{-1} \ell h)r \,=\, L\ell hr \,=\, Lhr.
	\]
	
\\ To see that $f$ is a morphism, note first that $f(Kr \rac m) = f(Kr) \rac m$ holds by associativity if both sides are defined, so it remains to check the definedness condition.  For this, we need to check that $rm \in R_e$ if and only if $hrm \in R_e$. But this follows from the fact that $H_e$ acts on the left of~$R_e$. Finally, note that $f(K) = f(Ke) = Lhe = Lh$, as desired.
\end{proof}

\begin{Cor}\label{conjcond}
$R_e/K$ is isomorphic to $R_e/L$ if and only if $K$ and $L$ are conjugate in $H_e$.
\end{Cor}
\begin{proof}
Suppose $R_e/K$ and $R_e/L$ are isomorphic. Applying the preceding proposition to the isomorphism and its inverse, we find that $K \setin x^{-1}Lx$ and $L \setin y^{-1}Ky$ for some $x,y \in H_e$. By finiteness, it follows that $K = x^{-1}y^{-1}Kyx$, and thus $x^{-1}Lx \setin x^{-1}y^{-1}Kyx = K$, so $K = x^{-1}Lx$ as desired.

	Conversely, if $K = x^{-1}Lx$, then $L = xKx^{-1}$, so by the proposition there exist morphisms $f$ and~$g$ with $f(K) = Lx$ and $g(L) = Kx^{-1}$. It is obvious that $f$ and~$g$ are inverses, hence $R_e/K$ and~$R_e/L$ are isomorphic.
\end{proof}

	For each $e \in \idrep$, fix a set $\subrep_e$ of representatives of the conjugacy classes of subgroups of~$H_e$, and let $K$ range over the subgroups in~$\subrep_e$. Every strong orbit determined by maximal subgroups is therefore isomorphic to precisely one quotient $R_e/K$. In particular, the classes $[R_e/K]$ freely generate~$\mathcal G$.
		
	If we take $M = H_e$ in Corollary~\ref{conjcond}, we obtain the standard result that $H_e/K$ and $H_e/L$ are isomorphic if and only if $K$ and $L$ are conjugate.  It follows that for each $e \in \idrep$, the classes $[H_e/K]$ with $K \in \subrep_e$ form a basis for the additive group of $\B(H_e)$.
	
	In this way, we see that $\G(M)$ and $\prod_e \B(H_e)$ have the same rank.
		
	If $M$ is distinguishable, then we saw that $\B(M) = \G(M)$, so in this case $\B(M)$ and $\prod_e \B(H_e)$ have the same rank.  On the other hand, if $M$ is not distinguishable, then $\B(M)$ has larger rank than~$\G(M)$: Indeed, if $u$ is any idempotent in a $\J$-class which is not distinguishable, then we saw above that $R_u/{\emax}$ is a strong orbit which is not isomorphic to any $R_e/K$.

\subsubsection{The structure theorem}

	Recall that $\phi : \B(M) \to \prod_e \B(H_e)$ given by $\phi([X]) = (e: [Xe])$ is a ring homomorphism. We now order the bases given in the previous section and show that the matrix for $\phi$ restricted to $\G(M)$ with respect to these ordered bases is upper-triangular with $1\rm s$ on the diagonal.
		
	Since $\idrep$ contains one idempotent per $\J$-class, the relation $\gej$ defines a partial order on~$\idrep$. Via topological sorting, refine this partial order into a linear order so that $e$ precedes~$f$ if $e \gneqq_{\J} f$. Extend the ordering of idempotents to an ordering of basis elements by arbitrarily choosing an order of the representative subgroups in each $\subrep_e$.
			
	To show that the matrix of $\phi$ restricted to $\G(M)$ with respect to these ordered bases is upper-triangular with $1\rm s$ on the diagonal, it suffices to show:
	
	\cl (i) $(R_e/K)e \,=\, H_e/K$ \\
	
\\ so that the $(e,e)\hbox{-block}$ of the matrix is an identity matrix, and for $f \in \idrep$:

	\cl (ii) $(R_e/K)f \,=\, \emptyset$ \q if $e$ precedes $f$ \\
	
\\ so that the matrix is $0$ below the diagonal.

	As for (i), the elements of $(R_e/K)e$ are $Kre$, where $re \in R_e$. By stability, $re \in H_e$, so $Kre \in H_e/K$, and thus $(R_e/K)e \setin H_e/K$. Conversely, if $h \in H_e$ then certainly $h \in R_e$, so $Kh = Khe$ is an element of $(R_e/K)e$. This shows that $(R_e/K)e = H_e/K$, as desired.
	
	As for (ii), if $e$ precedes $f$, then $f \not\gej e$ by construction, hence $f$ annihilates~$R_e$, and thus $(R_e/K)f = \emptyset$, which is what we wished to show.
	
	Thus $\phi$ restricted to $\G(M)$ is an isomorphism of groups, and we have established the following theorem:

\begin{Thm}
The map $\phi : \B(M) \to \prod_e \B(H_e)$ is a surjective ring homomorphism that restricts to an isomorphism on $\G(M)$, the subgroup of $\B(M)$ generated by the classes of strong orbits which are determined by maximal subgroups. \qed
\end{Thm}

\begin{Cor}
$M$ is distinguishable if and only if $\B(M) \cong \prod_e \B(H_e)$.
\end{Cor}
\begin{proof}
If $M$ is distinguishable, then $\B(M) = \G(M) \cong \prod_e \B(H_e)$ by Corollary~\ref{burnsideg}. If $M$ is not distinguishable, then, as discussed earlier, the ranks of $\B(M)$ and $\prod_e \B(H_e)$ differ, so these rings are not isomorphic.
\end{proof}

	We will see an example below where $M$ is not distinguishable, yet nevertheless the Burnside rings of the maximal subgroups determine much of the structure of the Burnside ring of~$M$.

\subsection{Examples}

	We now apply our structure theorem to a few monoids of note.

\subsubsection{$\ER$ monoids}

	A monoid $M$ is \term{$\ER$} if every $\R$-class of $M$ contains at most one idempotent.  Semigroups satisfying this property were studied by Rhodes and~Tilson in \cite[Definition~2.1]{rhodestilson}, and further explored by Stiffler in \cite[Theorem~3.18]{stiffler}.  (This property was called~``$R_1$'' in these early papers.)  More recent work can be found for example in \cite[Section~4.8]{qtheory} and \cite{msquiver}.

	Let $J$ be a regular $\J$-class of an $\ER$ monoid, and let $L$ and $L'$ be distinct $\L$-classes in~$J$. By a standard fact about finite monoids, every $\L$-class in a regular $\J$-class contains an idempotent.  Thus we may suppose that $L$~contains an idempotent~$e$, and hence $L \cap R_e$ contains~$e$. But $L'$ does not contain $e$, and $R_e$ does not contain any idempotents besides~$e$, so $L' \cap R_e$ does not contain an idempotent, and hence $R_e$ distinguishes $L$ and~$L'$. Thus we have shown that $\ER$ monoids are distinguishable, and hence our structure theorem applies to any $\ER$ monoid.
	
	The class of $\ER$ monoids contains other commonly-encountered classes of monoids which we mention here:
	
	We say that $M$ has \term{commuting idempotents} if $e \kern-.8pt f = f \kern-1.1pt e$ for all idempotents $e,f \in M$. Special cases include monoids with central idempotents and commutative monoids.  If $M$ has commuting idempotents and $e \R f$, then $e = f \kern-1.1pt e = e \kern-.8pt f = f$, so $M$ is $\ER$, and hence our structure theorem applies to monoids with commuting idempotents.
	
	We say that $M$ is a \term{block group} if for all $x \in M$, there is at most one $y \in M$ such that $xyx = x$ and $yxy = y$. If such a $y$ exists, we call $y$ the \term{inverse} of $x$. Suppose that $e$ and $f$ are idempotents in the same $\R$-class of a block group.  Then $e \kern-1pt f \kern-1pt e = e$ and $f \kern-1.1pt e \kern-.8pt f = f$, which shows that $e$ and $f$ are inverses.  But we also obviously have $eee = e$, which shows that $e$ is its own inverse.  Since inverses are unique in block groups, we must have $e = f$. Thus each $\R$-class contains at most one idempotent, which means that block groups are $\ER$, and hence our structure theorem applies to block groups.
	
	Finally, we say that a monoid is an \term{inverse monoid} if every element has a unique inverse.  Obviously inverse monoids are block groups, so our structure theorem applies to them as well.

\subsubsection{The monoid of square matrices over a finite field}

	In this section, we consider $M_{n,q}$, the monoid of $n \times n$ matrices over a finite field~$F$ with $|F| = q$ under matrix multiplication. We view $M_{n,q}$ as acting on the right of $F^n$, the space of $1 \times n$ row vectors over the same field.
	
	First we give easily-verified translations of Green's relations:

	\medbreak
	\ul \hbox to 4.5em{$X \gej Y$\hfil} $\iff$ \q $\rank(X) \ge \rank(Y)$\\
	\ul \hbox to 4.5em{$X \gel Y$\hfil} $\iff$ \q $\im(X) \contains \im(Y)$\\
	\ul \hbox to 4.5em{$X \ger Y$\hfil} $\iff$ \q $\ker(X) \setin \ker(Y).$\\
	\medbreak

	Now we show that $M_{n,q}$ is distinguishable.  Let $J$ be a regular $\J$-class, so that all matrices in~$J$ have the same rank, and let $X,Y \in J$ lie in distinct $\L$-classes, so that $X$ and~$Y$ have distinct images.  Let $v_1X,\dots,v_kX$ be a basis for~$\im(X)$, and let $wY$ be any vector in $\im(Y) \setminus \im(X)$. Then $wY$ and the vectors $v_iX$ are linearly independent, so there is a matrix $D \in J$ which satisfies $(v_iX)D = v_iX$ and $(wY)D = 0$. It is now clear that $\rank(XD) = k$, so $XD \in J$, but $\rank(YD) < k$, so $YD \notin J$. Hence by Clifford-Miller, $R_D$ distinguishes $L_X$ and~$L_Y$, and thus $M_{n,q}$ is distinguishable, as claimed.  It follows that our structure theorem applies to $M_{n,q}$.  It is easy to show that a maximal subgroup of rank~$r$ is isomorphic to~$GL_{r,q}$, so we have:
	
	\[
	\B(M_{n,q})
		\,\cong\, \prod_{r=0}^n \, \B(GL_{r,q}).
	\]

\subsubsection{The full transformation monoid}

	We now present a more complicated and interesting example.  Consider $T_n$, the set of total maps from $X = \{1,\dots,n\}$ to itself under the operation of function composition.  We view $T_n$ as acting on the right of~$X$.
	
	As we will see, $T_n$ has one $\J$-class which is not distinguishable, so our structure theorem does not directly apply to $T_n$. Nevertheless, we will be able to determine the structure of $\B(T_n)$ along similar lines.

\customhead{Green's relations on $T_n$}
	
	Define $\im(f) = \{ xf : x \in X \}$ and let $\rank(f)$ be the size of $\im(f)$. Further define an equivalence relation $\ker(f)$ by $(x,y) \in \ker(f)$ if and only if $xf = yf$. The equivalence classes of $\ker(f)$ are called the \term{fibers} of~$f$.
	
	The following translations are easily verified:
	
	\medbreak
	\ul \hbox to 4.2em{$f \gej g$\hfil} $\iff$ \q $\rank(f) \ge \rank(g)$\\
	\ul \hbox to 4.2em{$f \gel g$\hfil} $\iff$ \q $\im(f) \contains \im(g)$\\
	\ul \hbox to 4.2em{$f \ger g$\hfil} $\iff$ \q $\ker(f) \setin \ker(g).$\\
	\medbreak

	Write $J(r)$ for the $\J$-class of transformations of rank~$r$, $1 \le r \le n$. Each $\J$-class is regular, and we fix an idempotent from each $\J$-class. Let $e_r$ be defined by:
	
	\[
	xe_r \ = \ 
		\begin{cases}
			\,x & {\rm if} \ \ 1 \le x \le r \\
			\,r & {\rm if} \ \ r+1 \le x \le n \\
		\end{cases} 
	\]
	
\\ We may think of $e_r$ as an identity function, where $r$ through $n$ have been collapsed to the single number~$r$. It is trivial to check that $e_r$ is an idempotent of rank~$r$. Let $L(r) = L_{e_r}$ and similarly let $R(r) = R_{e_r}$ and $H(r) = H_{e_r}$. From our translation of Green's relations, it is clear that $H(r)$ is isomorphic to $\Sym(r)$, the symmetric group on $r$~letters.
	
\customhead{Distinguishability of $\J$-classes in $T_n$}

	We now turn to the question of which $\J$-classes in $T_n$ are distinguishable.
		
	If a function in $T_n$ has rank~$n$, then its image is~$X$, so $J(n)$ comprises a single $\L$-class and is hence trivially distinguishable.
	
	On the other hand, we show $J(1)$ is only distinguishable when $n=1$. From Green's relations, $J(1)$ consists of the constant functions.  Let $\kappa_i$ denote the constant function with value~$i$. For any $f \in T_n$ we have $f \kappa_i = \kappa_i$, so each constant function is a right zero-element and forms a distinct $\L$-class. But the product of any two constant functions is constant, and thus the product of any two elements of $J(1)$ lies in~$J(1)$. Hence by Clifford-Miller, all $\L$-classes in~$J(1)$ are indistinguishable, and so $J(1)$ is not distinguishable if $n \ge 2$.
		
	Finally, we show that $J(r)$ is distinguishable if $r > 1$. Let $f$ and $g$ be functions in $J(r)$ which lie in distinct $\L$-classes of~$J(r)$, so that $f$ and $g$ have distinct images.  Choose $c_1,\dots,c_r$ so that $\im(f) = \{c_1,\dots,c_r\}$ but $c_1 \notin \im(g)$, and then choose $c_{r+1},\dots,c_n$ so that $X = \{c_1,\dots,c_n\}$. Let the function $k$ be defined by:

	\[
	c_ik \ = \  
		\begin{cases}
		\,c_i & {\rm if} \ \ 1 \le i \le r \\
		\,c_r & {\rm if} \ \ r+1 \le i \le n \\
		\end{cases} 
	\]
	
\\ By construction, $k$ has rank $r$ and so does $fk = f$, so both $k$ and $fk$ are in~$J(r)$. On the other hand, we have $\im(gk) \setin \{c_2,\dots,c_r\}$, so $\rank(gk) < r$ and hence $gk \notin J(r)$. By Clifford-Miller, $J(r)$ is distinguishable.
	
	Since $J(1)$ is not distinguishable for $n \ge 2$, our structure theorem does not directly apply to~$T_n$, except in the trivial case $n = 1$.

\customhead{The strong orbits of $T_n$}

	In this section, we show that $T_n$ has only one strong orbit class which is not determined by maximal subgroups: namely, the identity of $\B(T_n)$.
		
	A strong orbit of $T_n$ is isomorphic to a quotient of $R(r)$, where $1 \le r \le n$. For $r > 1$, $J(r)$ is distinguishable, and thus Corollary~\ref{distcong} shows that such quotients are determined by maximal subgroups.
	
	Now consider quotients of $R(1)$. Since constant functions are right zero-elements, $R(1)$~contains all constant functions.   If the functions $\kappa_a$ and~$\kappa_b$ are distinct and congruent, and $\kappa_c$ is an arbitrary constant function, then right-multiplying by a function in~$T_n$ which fixes~$a$ and sends $b$ to~$c$ shows that $\kappa_a$ and $\kappa_c$ are congruent as well.  Thus the only congruences on $R(1)$ are $=$ and $\emax$, and hence the only quotients of~$R(1)$ are $R(1)$ itself and $R(1)/{\emax}$.
	
	In fact, $R(1)/{\emax}$ is isomorphic to $\bf 1$, the identity of~$\B(T_n)$, which recall consists of a single point fixed by all of~$T_n$. Indeed, we saw that each constant function forms its own $\L$-class and that these $\L$-classes are indistinguishable, so $R(1)/{\emax}$ contains a single congruence class consisting of all constant functions. If $f$~is any function in~$T_n$, we have $\kappa_a f = \kappa_{af}$, which shows that the action of $T_n$ fixes this class, and hence $R(1)/{\emax}$ is isomorphic to $\bf 1$.
	
	Thus we have shown that the identity of~$\B(T_n)$ is the only strong orbit class which is not determined by maximal subgroups.

\customhead{The Burnside ring of $T_n$}

	To compute the Burnside ring of $T_n$, we need a lemma.

\begin{Lemma}
For any monoid~$M$, the additive subgroup of~$\B(M)$ generated by the nonidentity strong orbit classes forms an ideal.
\end{Lemma}
\begin{proof}
It suffices to show that the product of a nonidentity strong orbit and any strong orbit does not contain a strong orbit isomorphic to~$\bf 1$; in other words, that such a product does not contain any points fixed by~$M$. Indeed, let $X$~be a nonidentity strong orbit and let $Y$ be any strong orbit.  Then for all $x \in X$ there exists $m \in M$ such that either $xm$ is undefined or $xm$ is defined and $xm \ne x$. Hence for all $y \in Y$ either $(x,y)m$ is undefined or $(x,y)m = (xm,ym) \ne (x,y)$. Therefore no point of~$X \sprod Y$ is fixed by all of~$M$, which is what we wished to show.
\end{proof}

	Now let $I$ be the additive subgroup of $\B(T_n)$ generated by the nonidentity strong orbit classes.  As we demonstrated above, these are the classes of strong orbits determined by maximal subgroups; in other words, we have $I = \G(T_n)$. By the lemma, $I$ is an ideal of $\B(T_n)$. Thus the map $\phi$ from our structure theorem, restricted to~$I = \G(T_n)$, is a ring isomorphism.  It follows that $I \cong \prod_r \B(H(r))$ is a unital ring, and hence:
	
	\[
	\B(T_n)
	\, \cong \, \B(T_n)/I \times I
	\, \cong \, \B(T_n)/I \times \prod_{r=1}^n \B(H(r)).
	\]
	
\\ Finally, observe that the map from $\B(T_n)$ to $\Z$ which sends an element of $\B(T_n)$ to its coefficient of $\bf 1$ is a surjective ring homomorphism whose kernel is~$I$. Hence $\B(T_n)/I \cong \Z$. Remembering that $H(r) \cong \Sym(r)$, we have shown at last:
	
	\[
	\B(T_n)
		\, \cong \, \Z \times \prod_{r=1}^n \B(\Sym(r)).
	\]

\section{The Burnside algebra and the table of marks}

	In this section, we discuss some well-known concepts from representation theory, namely the \term{Burnside algebra} over the field $\Q$ of rational numbers, introduced by Solomon in~\cite{solomon}, as well as what Burnside in~\cite{burnside} calls the \term{table of marks}. In order to translate this theory from the language of groups to that of monoids, we will need the notion of a \term{lax morphism}, to be introduced in due course.

\subsection{Burnside algebra}

	Let $\bk$ be a field.  The \term{Burnside algebra} of $M$ over~$\bk$ is defined to be the tensor product $\Bk = \bk \otimes_\Z \B$. For $k \in \bk$ and $[X] \in \B$, we write $k[X]$ instead of $k \otimes [X]$, so that an arbitrary element of~$\Bk$ may be written $\sum k_i[X_i]$. Because the strong orbit classes form a basis for the additive group of $\B$, they form a basis for $\Bk$ over $\bk$ as well.  Hence for a finite monoid, $\Bk$ is a commutative finite-dimensional $\bk$-algebra.

	We now restrict our attention to the case $\bk = \Q$. It is known that the Burnside algebra of a group over $\Q$ is isomorphic to the direct product of copies of $\Q$; in other words, $\BQ(G)$ is \term{semisimple}. Our goal in this section is to show that the same is true for monoids in general.
	
\subsection{Marks}

	In the group setting, one can prove that $\BQ(G)$ is semisimple via a mapping known as the \term{table of marks}. We briefly summarize that theory: For a subgroup $H$ of $G$ and a $G$-set~$X$, we define $X^H = \{ x \in X : xh = x \ \hbox{for\ all}\ h \in H \}$ to be the subset of elements of~$X$ which are fixed by every element of~$H$. The number $|X^H|$ is known as the \term{mark} of~$H$ on~$X$.
	
	The reason why this definition is meaningful for groups, but not monoids, will become apparent in a moment.  Let $X$ be any $G$-set, and let $\Omega$ be a transitive $G$-set with $\a \in \Omega$. A morphism $f : \Omega \to X$ is completely determined by where it sends $\a$, since $f(\a g) = f(\a) g$. For $f$ to be well-defined, it is necessary that:
	
	\[
	\a g = \a \qq{\rm implies}\qq f(\a)g = f(\a g) = f(\a) ;
	\]
	
\\ in other words, if $g$ fixes $\a$ then $g$ fixes $f(\a)$. This condition is also sufficient for $f$ to be well-defined: Indeed, if the condition holds, then $\a g = \a h$ implies that $gh^{-1}$ fixes $\a$, so $gh^{-1}$ fixes $f(\a)$, and hence $f(\a)g = f(\a)h$. But then $f(\a g) = f(\a h)$, which shows that $f$ is well-defined.

	Now take $\Omega = G/H$ and $\a = H$ in the above.  The elements that fix $\a$ are precisely the elements of $H$, hence $f$ is well-defined if and only if $f(\a)$ is fixed by $H$; in other words, if and only if $f(\a) \in X^H$. Thus we see that the number of morphisms from $G/H$ to $X$ equals the mark of $H$ on $X$.
	
	If we replace the group with a monoid, the condition we identified above is necessary but no longer sufficient to guarantee that $f$ is well-defined.  Indeed, following the above argument, if $\a m = \a n$, we would like to conclude that either $mn^{-1}$ or $nm^{-1}$ fixes $\a$, but it may be that neither~$m$ nor~$n$ are invertible.  The definition of mark in terms of fixed points is therefore a convenient phrasing which is suitable for groups, but not for monoids in general.

\subsection{Lax morphisms}

	This discussion suggests we should calculate marks in the monoid setting by counting the number of functions defined on a strong orbit $\Omega$ which are completely determined by their behavior on one element $\a \in \Omega$. In symbols, this means that wherever $\a \cdot m$ is defined, then $f(\a) m$ is also defined and given by $f(\a \cdot m) = f(\a) m$.
	
	More generally, we say that a map $f : X \to Y$ of $M$-sets is a \term{lax morphism} if $f(xm) = f(x)m$ holds wherever $xm$ is defined. Note that we do not require both sides to be undefined when $xm$ is undefined, as would be required if $f$ were a morphism.

	Let $\Lax(X,Y)$ be the set of lax morphisms from $X$ to~$Y$, and note that, in the group setting, the mark of a subgroup $H$ of $G$ on a $G$-set~$X$ is $|\Lax(G/H,X)|$.

	Every morphism is a lax morphism, so in particular isomorphisms are lax morphisms.  It is also easy to verify that the composition of lax morphisms is a lax morphism.  Thus if $X \cong X'$ and $Y \cong Y'$, $\Lax(X,Y)$ is in bijective correspondence with $\Lax(X',Y')$ by pre- and post-composition with the appropriate isomorphism.
	
	To show that a lax morphism $f$ is a morphism, it has to be checked that $xm$ is defined if $f(x)m$ is defined, which is not true in general.  Indeed, we saw in earlier sections that inclusion of a subquotient and projection are not generally morphisms.  Let us revisit those examples now.
	
	If $S$ is a subquotient of an $M$-set $X$, then the inclusion map $\iota : S \to X$ is a morphism if and only if $S$ is invariant. However, $\iota$ is always a lax morphism, since if $s \cdot m$ is defined then we have $s \cdot m = sm$, and hence $\iota(s \cdot m) = \iota(sm) = sm = \iota(s)m$ for all $s \in S$ and $m \in M$.
	
	The projection maps:
	
	\medbreak
	\ul \hbox to 8.0em{$\pi_X : X \sprod Y \to X$,\hfil} $\pi_X((x,y)) \,=\, x$\\
	\ul \hbox to 8.0em{$\pi_X : X \sprod Y \to Y$,\hfil} $\pi_X((x,y)) \,=\, y$\\
	\medbreak

\\ are also lax morphisms: Indeed, if $(x,y)m$ is defined, then $xm$ and $ym$ are both defined, so:

	\[
	\pi_X((x,y)m) \,=\, \pi_X((xm,ym)) \,=\, xm \,=\, \pi_X((x,y))m 
	\]
	
\\ and similarly for $\pi_Y$.

	Finally, observe that if $f$ is a lax morphism, then $x \ge x'$ implies $f(x) \ge f(x')$, so $x \sim x'$ implies $f(x) \sim f(x')$ as well.  Hence strong orbits map into, though not necessarily onto, strong orbits.

\subsection{Table of marks}

	Let $\cO$ range over a set of representatives for the strong orbit classes of $M$. Define a map~$\phi$ from strong orbit classes into $\prod_{\cO} \Z$ by:
	
	\cl $(*)$ $\phi([\Omega]) \,=\, (\cO: |\Lax(\cO,\Omega)|\, )$ \\
	
\\ where the notation $(\cO: P(\cO))$ stands for the tuple with $P(\cO)$ in the coordinate corresponding to~$\cO$. By the comments in Section~5.3, this map is well-defined and independent of the representatives~$\cO$ and $\Omega$.

	Since the strong orbit classes form a basis for the additive group of $\B$, we may extend this map uniquely to a $\Z\hbox{-linear}$ map $\phi : \B \to \prod_{\cO} \Z$, and in fact property $(*)$ holds for any $M$-set $X$:  Indeed, if $\Omega_1,\dots,\Omega_n$ are the strong orbits of~$X$, then $[X] = [\Omega_1]+\cdots+[\Omega_n]$ in $\B$, and hence:
	
	\[
	\phi([X]) \,=\, (\cO : |\Lax(\cO,\Omega_1)|+\cdots+|\Lax(\cO,\Omega_n)|\,).
	\]
	
\\ But it is clear that:
	
	\[
	|\Lax(\cO,\Omega_1)|\,+\,\cdots\,+\,|\Lax(\cO,\Omega_n)| \,=\, |\Lax(\cO,X)|
	\]
	
\\ since a lax morphism from $\cO$ to $\Omega_i$ followed by inclusion into~$X$ is a lax morphism from $\cO$ to $X$; and conversely, every lax morphism from $\cO$ to $X$ is a lax morphism from $\cO$ to some $\Omega_i$ by our observation above that lax morphisms map into strong orbits. These correspondences are clearly inverses, and thus we have $\phi([X]) = (\cO : |\Lax(\cO,X)|\,)$, as claimed.

	We will show that $\phi$ is an injective $\Z\hbox{-algebra}$ homomorphism by ordering the basis $\{[\cO]\}$ and showing that the matrix of $\phi$ with respect to this ordered basis is upper-triangular with nonzero entries on the diagonal.  This matrix is called the \term{table of marks}.

\subsection{A partial order on strong orbit classes}

	In the group setting, strong orbits are isomorphic to systems of cosets of subgroups.  Thus to order the basis, one may choose a representative subset from each conjugacy class, and define $G/H \succeq G/K$ if and only if $H \le g^{-1}Kg$ for some $g \in G$. Clearly this is unsuitable for our purposes.  However, by Proposition~\ref{morphcond} in Section~4.3.2, we know that this is equivalent to saying that there exists a morphism from~$G/H$ to~$G/K$; or, in our new terminology, that $|\Lax(G/H,G/K)| \ne 0$.
	
	Let us therefore define an ordering $\succeq$ on strong orbits by $\cS \succeq \cT$ iff $|\Lax(\cS,\cT)| \ne 0$. Because the identity map is a lax morphism and the composition of lax morphisms is a lax morphism, $\succeq$ is reflexive and transitive, hence a preorder.  Thus if we define $\cS \cong^* \cT$ if and only if $\cS \succeq \cT$ and $\cT \succeq \cS$, then $\cong^*$ is an equivalence relation, and $\succeq$ induces a partial order on the equivalence classes of $\cong^*$.
	
	We will show that the equivalence classes of $\cong^*$ are just the usual isomorphism classes of strong orbits, and hence that $\succeq$ gives a partial order on strong orbit classes.

\begin{Lemma}\label{laxaut}
A lax endomorphism of a strong orbit is an automorphism.
\end{Lemma}

\begin{proof}
By Proposition \ref{apex}, let $\Omega = \a R_e$ be a strong orbit, where $e$ is an idempotent in the apex such that $\a e = \a$ and $\a r$ is defined for all $r \in R_e$. Let $f : \Omega \to \Omega$ be a lax morphism.  Since $e$ fixes $\a$, we have:

	\[
	f(\a) \,=\, f(\a \cdot e) \,=\, f(\a) \cdot e
	\]
	
\\ so $e$~fixes $f(\a)$. Let $f(\a) \cdot m$ be an arbitrary element of $\Omega$. Then $f(\a) \cdot m = f(\a) \cdot em$, so $em$ does not annihilate $\Omega$ and hence $em \gej e$. We always have $e \gej em$, so $e \J em$ and hence $em \in R_e$ by stability.  Thus $\a \cdot em$ is defined, and so:

	\[
	f(\a \cdot em) \,=\, f(\a) \cdot em \,=\, f(\a) \cdot m 
	\]
	
\\ which shows that $f$ is a surjection.  By finiteness, $f$ is a bijection.

	To show that $f$ is a morphism, we must show that $f(\a r)m \in \Omega$ implies $\a rm \in \Omega$ for all $r \in R_e$. Now $f(\a r)m = f(\a)rm$, so $rm$ does not annihilate~$\Omega$, hence $rm \gej r$. But we always have $r \gej rm$, so $r \J rm$ and hence $r \R rm$ by stability.  Thus $rm \in R_e$ and hence $\a rm \in \Omega$, as desired.
\end{proof}

\begin{Thm}\label{isoiso}
For strong orbits $\cS$ and $\cT$, we have $\cS \cong^* \cT$ if and only if $\cS \cong \cT$.
\end{Thm}

\begin{proof}
First suppose that $f : \cS \to \cT$ is an isomorphism of strong orbits. Then $f$ and $f^{-1}$ are lax morphisms, so $\cS \succeq \cT$ and $\cT \succeq \cS$, hence $\cS \cong^* \cT$.
	
	Conversely, suppose $\cS \cong^* \cT$ so that $\cS \succeq \cT$ and $\cT \succeq \cS$. Then there exist lax morphisms $f : \cS \to \cT$ and $g : \cT \to \cS$. By the above lemma, $gf : \cS \to \cS$ and $fg : \cT \to \cT$ are automorphisms, so $f$ is a bijection by set theory.  To see that $f$ is a morphism, we have to show that $\cs \cdot m$ is defined if $f(\cs) \cdot m$ is, so suppose $f(\cs) \cdot m$ is defined for $\cs \in \cS$. Then because $g$ is a lax morphism and $gf$ is a morphism, we have:
	
	\[
	g(f(\cs) \cdot m) \,=\, g(f(\cs)) \cdot m \,=\, (gf)(\cs) \cdot m \,=\, (gf)(\cs \cdot m) 
	\]
	
\\ hence $\cs \cdot m$ is defined.  This shows that $f$ is a morphism, hence an isomorphism, so $\cS \cong \cT$ as desired.
\end{proof}

	As we stated above, $\succeq$ gives a partial order on the equivalence classes of~$\cong^*$. But by Theorem~\ref{isoiso}, these are just isomorphism classes of strong orbits, hence $\succeq$ gives a partial order on strong orbit classes, as desired.

\subsection{$\BQ$ is semisimple}

	Recall that $\phi : \B \to \prod_{\cO} \Z$ is a $\Z\hbox{-linear}$ map which satisfies:
	
	\[
	\phi([X]) \,=\, (\cO: |\Lax(\cO,X)|\, ) 
	\]
	
\\ for all pointed $M$-sets $X$.

	We show that $\phi$ is an injective $\Z\hbox{-algebra}$ homomorphism. The only lax morphism from a strong orbit~$\cO$ to~$\bf 1$ is the constant map, hence:
	 
 	\[
	\phi({\bf1}) \,=\, (\cO: |\Lax(\cO,{\bf 1})|\, ) \,=\,
		 (\cO: 1) 
	\]
	
\\ which is the identity element of $\prod_{\cO} \Z$, so $\phi$ preserves the identity.  To see that $\phi$ is multiplicative, it suffices to show:
	 
	\[
	|\Lax(\cO,X \sprod Y)| \, = \, |\Lax(\cO,X)| \,\cdot\, |\Lax(\cO,Y)|
	\]
	 
\\ for all strong orbits $\cO$ and $M$-sets $X$ and $Y$. Recall that projection maps are lax morphisms, so if $f : \cO \to X \sprod Y$ is a lax morphism, then so are $\pi_X \circ f : \cO \to X$ and $\pi_Y \circ f : \cO \to Y$. Conversely, if $g : \cO \to X$ and $h : \cO \to Y$ are lax morphisms, it is straightforward to check that $(g,h) : \cO \to X \sprod Y$ given by $(g,h)(\a) \,=\, (g(\a),h(\a))$ is a lax morphism.  Moreover, these correspondences are inverses of each other, in the sense that:
	
	\[
	f = (\pi_X \circ f,\pi_Y \circ f) \qq{\rm and}\qq 
	\begin{cases}
	\,g \,=\, \pi_X \circ (g,h)\\
	\,h \,=\, \pi_Y \circ (g,h)
	\end{cases}.
	\]
	
\\ Thus $|\Lax(\cO,X \sprod Y)| = |\Lax(\cO,X)| \cdot |\Lax(\cO,Y)|$, which shows that $\phi$ is multiplicative, and hence a $\Z\hbox{-algebra}$ homomorphism.

	Now we order our basis and construct the table of marks.  By topological sorting, refine the partial order $\succeq$ on the set $\{\cO\}$ of strong orbit representatives into a linear order, so that $\cO$~precedes $\cO'$ if $\cO \succ \cO'$. Write the representatives as $\{\cO_i\}$, where $i$ ranges from $1$ to the rank of $\B$, so that $\{[\cO_i]\}$ is an ordered basis for~$\B$.
	
	Consider the matrix for $\phi$ with respect to this ordered basis.  The $i,j\hbox{-entry}$ of this matrix is $|\Lax(\cO_i,\cO_j)|$. If $i > j$, then $\cO_i \nsucceq \cO_j$, hence $|\Lax(\cO_i,\cO_j)| = 0$. This shows that the matrix for $\phi$ is upper-triangular.  Hence we have:
	
\begin{Thm}
The map $\phi : \B \to \prod_\cO \Z$~is an injective $\Z$-algebra homomorphism, and hence $\B$ is isomorphic to a subring of $\prod_\cO \Z$ with index $\prod_\cO |\!\Aut(\cO)|$ as a subgroup.
\end{Thm}
\begin{proof}
We showed that $\phi$ is a $\Z$-algebra homomorphism, and that the matrix for $\phi$ is upper-triangular.  The diagonal entries of this matrix are $|\Lax(\cO_i,\cO_i)|= |\!\Aut(\cO_i)| \ne 0$, so $\phi$ is injective, and the index of the image is the product of the diagonal entries, namely $\prod_\cO |\!\Aut(\cO)|$.
\end{proof}

\begin{Cor}
$\BQ$ is semisimple.
\end{Cor}
\begin{proof}
The matrix of $\phi$ is upper-triangular with nonzero entries on the diagonal, hence it is invertible over $\Q$. It follows that $\phi : \BQ \to \prod_i \Q$ is a $\Q\hbox{-algebra}$ isomorphism.  Since $\BQ$ is commutative and finite-dimensional, $\BQ$ is semisimple as claimed.
\end{proof}

\section*{Acknowledgements}

	The research which led to this paper was funded by the Barnett and Jean Hollander Rich Summer Mathematics Internship, which I received in Summer 2022 and~2023 while I was a Master's student at the City College of New York.  I thank the City College mathematics department sincerely for their financial support.

	My mentor for these summers was Benjamin Steinberg.  Professor Steinberg suggested exploring a Burnside ring for monoids by analogy with partial transformation monoids.  He introduced me to the notions of monoid actions, $M$-sets, as well as the fundamentals of semigroup theory.  He pointed me to books and articles, and suggested many avenues for exploration and ideas for proofs.  While the work contained in these pages is my own, I am deeply indebted to Professor Steinberg for giving life and shape to this project, and for the generous way he has always been willing to explain a concept, critique a proof, or provide detailed feedback on my writing.

\appendix
\setcounter{secnumdepth}{0}
\section{Appendix}

	Let $R_e$ be the $\R$-class of an idempotent~$e$, and let $\equiv$ be a right congruence on~$R_e$. In Section 3, we defined $\fL$ to be the subgroup of $H_e$ consisting of all elements of $H_e$ which preserve~$\equiv$ by left multiplication.  Thus $\fL$ acts on $R_e/{\equiv}$ by automorphisms.  Letting $\fK$ be the kernel of this action, we showed that $\Aut(R_e/{\equiv}) \cong \fL/\fK$. By analogy with a similar result from group theory, it is reasonable to conjecture that $\fL = \mathcal{N}_{H_e}(\fK)$. However, Benjamin Steinberg constructed a counterexample, which we present here.
	
	\newcommand{\mat}[9]{%
		\left(\!
			\begin{smallmatrix}
				#1 & #2 & #3\\#4 & #5 & #6\\#7 & #8 & #9
			\end{smallmatrix}
		\!\right)
	}

	Let $M = 
	\left\{
		\mat100010001\right.$, $\pm\mat100100100$, $\pm\mat010010010$,
		$\pm\mat001001001$, $\pm\mat{\phantom-1}00{-1}00{-1}00$,
		$\pm\mat0{\phantom-1}00{-1}00{-1}0$, $\left.\pm\mat00{\phantom-1}00{-1}00{-1}
	\right\}
	$ and $e = \mat100100100$. By direct computation, we have:

	\medbreak
	\ul \hbox to 1.2em{$J_e$} $=\ M \setminus \left\{\mat100010001\right\}$\\
	\ul \hbox to 1.2em{$R_e$} $=\ \left\{ \pm\mat100100100, \pm\mat010010010, \pm\mat001001001 \right\}$\\
	\ul \hbox to 1.2em{$H_e$} $=\ \left\{ \pm\mat100100100 \right\} \, \cong \, \Z/2\Z$. \\
	\medbreak
	
	Define an equivalence relation $\equiv$ on $R_e$ by $\mat010010010 \equiv \mat001001001$, and all other matrices are equivalent only to themselves.  This is a right congruence because the second and third rows of any matrix in~$M$ are the same.  Since:
	
	\[
	-\mat100100100\!\mat010010010 \, = \,
	-\mat010010010 \, \not\equiv \, 
	-\mat001001001 \, = \,
	-\mat100100100\!\mat001001001 
	\]
	
\\ we have $\fL = \{e\}$, while $\fK = \{e\}$ so that $\mathcal{N}_{H_e}(\fK) = H_e$. Hence $\fL \ne \mathcal{N}_{H_e}(\fK)$ in general.  Moreover, $\Aut(R_e/{\equiv}) \cong \fL/\fK$ is trivial, while $\mathcal{N}_{H_e}(\fK)/\fK \cong \Z/2\Z$.

\end{document}